\documentclass{article}
\usepackage[utf8]{inputenc}
\usepackage[left=3cm,right=3cm,top=2.25cm,bottom=2.25cm]{geometry} 
\usepackage{epsfig,graphicx}
\usepackage{xcolor}
\usepackage{amsmath,amsfonts,mathrsfs,amsfonts,ulem, bbold,bm,bbm,blindtext}
\usepackage{algorithmicx,algorithm,algpseudocode}
\algnewcommand{\LineComment}[1]{\State \textcolor{blue}{\(\triangleright\) #1}}
\algnewcommand{\LineCommentg}[1]{\State \textcolor{gray}{\(\triangleright\) #1}}

\usepackage{xr-hyper}
\usepackage{hyperref}

\usepackage{multirow}

\renewcommand{\L}{\textrm{L}}

\newcommand{\M}{\textrm{M}}

\newcommand{\Prol}{\textrm{Prol}}
\newcommand{\Circ}{\textrm{Circ}}

\newcommand{\kprol}{k_{\textrm{prol}}}
\newcommand{\ktr}{k_\textrm{tr}}
\newcommand{\Edrug}{E_\textrm{drug}}

\newcommand{\Stem}{\textrm{Stem}}
\newcommand{\Cent}{\textrm{Cent}}

\newcommand{\ftr}{\textrm{ftr}}

\newcommand{\thetaTV}{\theta^\textrm{TV}}
\newcommand{\hatthetaTV}{\widehat{\thetaTV}}
\newcommand{\hatOmega}{\hat{\Omega}}

\newcommand{\thetaTVl}{\theta^{\textrm{TV}{(l)}}}

\newcommand{\rd}{\mathrm{d}}
\newcommand{\iid}{\mathrm{iid}}

\makeatletter
\newcommand{\objref}[4]{\def\obj@rg{#4}%
  #1\ifx\obj@rg\empty#2\else#3\xspace\ref{#4}--\fi\ref}
\makeatother
\newcommand{\Sobjref}[1]{\objref{#1}{~}{s}}

\usepackage{rotating} 
\newcommand{\Figref}[1][]{\Sobjref{Figure}{#1}}

\newcommand{\Eqref}[1][]{\Sobjref{Eq.}{#1}}

\newcommand{\Tabref}[1][]{\Sobjref{Table}{#1}}

\usepackage{adjustbox} 

\title{A continued learning approach for model-informed precision dosing: updating models in clinical practice}
\author{Corinna Maier$^{1,2}$, Jana de Wiljes$^{1}$,\\[1ex] Niklas Hartung$^{1}$, Charlotte Kloft$^{3}$, and Wilhelm Huisinga$^{1,\ast}$}

\begin{document}

\maketitle


\noindent
$^1$Institute of Mathematics, University of Potsdam, Germany\\[1ex]
$^2$Graduate Research Training Program PharMetrX: Pharmacometrics \& Computational Disease Modelling, Freie Universit\"at Berlin and University of Potsdam, Germany\\[1ex]
$^3$Department of Clinical Pharmacy and Biochemistry, Institute of Pharmacy, Freie Universit\"at Berlin, Germany\\
$^\ast$corresponding author\\
Institute of Mathematics, University of Potsdam\\
Karl-Liebknecht-Str. 24-25, 14476 Potsdam/Golm, Germany\\
Tel.: +49-977-59 33, Email: huisinga@uni-potsdam.de

\subsubsection*{Keywords}
Parameter estimation, Bayesian learning, inverse problems, Data Assimilation, Oncology, Chemotherapy, Therapeutic Drug Monitoring, Individualization, Personalized Medicine, Sequential Design, Target Concentration Intervention

\begin{abstract} 
Model-informed precision dosing (MIPD) is a quantitative dosing framework that combines prior knowledge on the drug-disease-patient system with patient data from therapeutic drug/ biomarker monitoring (TDM) to support individualized dosing in ongoing treatment.
Structural models and prior parameter distributions used in MIPD approaches typically build on prior clinical trials that involve only a limited number of patients selected according to some exclusion/inclusion criteria. 
Compared to the prior clinical trial population, the patient population in clinical practice can be expected to include also altered behavior and/or increased interindividual variability, the extent of which, however, is typically unknown. 
Here, we address the question of how to adapt and refine models on the level of the model parameters to better reflect this real-world diversity.
We propose an approach for continued learning across patients during MIPD using a sequential hierarchical Bayesian framework. The approach builds on two stages to separate the update of the individual patient parameters from updating the population parameters. Consequently, it enables continued learning across hospitals or study centers, since only summary patient data (on the level of model parameters) need to be shared, but no individual TDM data. 
We illustrate this continued learning approach with neutrophil-guided dosing of paclitaxel.
The present study constitutes an important step towards building confidence in MIPD and eventually establishing MIPD increasingly in everyday therapeutic use.
\end{abstract}


\section{Introduction}
Model-informed precision dosing (MIPD) is a quantitative framework for dose individualization based on modeling and simulation of exposure-response relationships integrating patient-specific data \cite{Keizer2018,Peck2021,Kluwe2020}.
The underlying models are developed based on clinical trial data typically using a nonlinear mixed effects (NLME) framework to describe the pharmacokinetics (PK) or pharmacodynamics (PD) of the drug and the variability between patients \cite{Lavielle2014}.
These PK/PD models allow to forecast important aspects of the therapy outcome based on patient characteristics ({\it a priori} predictions).
Therapeutic drug/biomarker monitoring (TDM) allows to further individualize model predictions ({\it a posteriori} predictions) and subsequently to adjust dosing.

When PK/PD models are used for MIPD, a `perfect model scenario' is generally assumed, in which the model represents the drug-patient-disease system sufficiently well, the variability of the outcome is adequately described and the prior study population (used to {\it develop} the model) is representative of the target individual patient (to which the model will be {\it applied}).
A certain model misspecification or population shift can, however, be expected due to the limited amount of data the models were built on. Specifically, the data from clinical trials only involve a limited number of patients, selected according to strict inclusion/exclusion criteria within a restricted time frame \cite{Keizer2018,Polasek2018,Zhao2013}. Therefore, models underlying MIPD will inevitably be confronted with deviating data in clinical routine, such as differences related to pathophysiology \cite{terHeine2020} or the patient population (comorbidities, comedications or special characteristics, e.g., morbidly obese, pregnant, or rare genotypes) \cite{Keizer2018,Polasek2018,Sheiner1992,Deitchman2019,Powell2020}.
In this `imperfect model scenario', the benefits of MIPD approaches may not reach their full potential.
It is therefore prudent to also improve the associated models as clinical routine data on the observed patient population becomes available.

For a given drug-disease-patient-system, there are often numerous models available within the literature, often based on different patient populations, e.g., for warfarin \cite{Hamberg2007,Ohara2014}, vancomycin \cite{terHeine2020,Uster2020} or ciclosporin \cite{Mao2018}.
Also, adjustments to the model used in a MIPD framework were necessitated after treatment of the first patient cohort \cite{Jodrell1994} or in retrospect \cite{Conley1989,Henrich2017}.
As an example of high clinical relevance, we focus on paclitaxel causing neutropenia as the most frequent and life-threatening toxicity in oncology.
Models describing paclitaxel-induced neutropenia build the basis for neutrophil-guided MIPD to individualize chemotherapy dosing \cite{Wallin2009a,Netterberg2017,Maier2020,Maier2020DARL}.
Since the publication of the gold-standard model for neutropenia \cite{Friberg2002}, many model variants have been developed, which differ not only in parameter estimates \cite{Kloft2006,Hansson2010,Joerger2007,Joerger2012}, but also in their structure \cite{Mensing2014,Pujo2016,Craig2017,Henrich2017}.

The challenge to choose between competing models is often approached via model averaging or model selection \cite{Uster2020}.
In model averaging, all candidate models are used, weighting the model predictions with the TDM data.
In contrast, in model selection, a single model is selected based on a retrospective external evaluation on independent data collected previously in the intended setting (from the same hospital and patient population) \cite{Zhao2013} 
and prospective fit-for-purpose verification \cite{terHeine2020}.
None of the approaches, however, integrates the new data collected during routine application of MIPD into the initial models underlying MIPD. 
Yet continued learning approaches based on an ever-growing amount of data have enormous potential to improve the predictive capabilities of MIPD in clinical practice.

The problem of transferability is a well-known problem in the machine learning literature; often called lifelong learning, continual learning \cite{Chen2015,Silver2013}, transfer learning \cite{PanYang2010,TorreyShavlik2010} or domain adaptation \cite{Jiang2008}.
Contrary to typical ML applications with direct access to big data, sensitive patient data may not be accessible and available to this extent across different sources.
Consequently current approaches based on pooling of data are not an option in this case, since they require direct access to TDM data of all patients \cite{Hughes2020}.
Hence, there is a need for approaches for continued model updating that are based on summary information of the data that is extracted locally and can be shared.


In this article, we propose an approach building on a sequential hierarchical Bayesian framework \cite{Lunn2013,Hooten2019} for continued model learning.
The underlying prior model used within MIPD is improved as new data from the target patient population are assimilated.
Importantly, the approach separates inference of the individual model parameters during a patient's therapy from the update of population parameters across patients.
The proposed approach is based on ideas from Bayesian integration of meta-analyses \cite{Lunn2013,Hooten2019}.
First, we demonstrate how a mismatch between the model and data-generating process could affect MIPD in an {\it in silico} trial setting in terms of model parameters and structural misspecifications.
We focus in the present article on how to adapt a model on the level of the typical model parameters and the magnitude of the inter-individual variability.
The proposed approach aims at bridging the gap between population analyses in model development and application of MIPD in therapeutic use. 

\section{METHODS}
\subsection{MIPD framework}
\label{sec:MIPDframework}

In MIPD a Bayesian framework is used. Here, we very briefly summarise the approach, for details see \cite{Maier2020,Maier2020DARL}.
Inference for a given ($i$-th) patient is based on a prior distribution $p(\cdot|\hatthetaTV,\hatOmega)$ and the likelihood $L(\cdot|y_i) = p(y_i|\cdot)$ of TDM data $y_i=(y_{i1},\dots,y_{1n_i})^T$, resulting in the posterior  $p(\cdot |y_i,\hatthetaTV,\hatOmega)$ of the individual parameters $\theta_i$ with
\begin{equation}\label{eq:Bayes}
p(\theta_i|y_i,\hatthetaTV,\hatOmega) \propto p(y_i|\theta_i)p(\theta_i|\hatthetaTV,\hatOmega)\,.
\end{equation}
The prior distribution in \Eqref{eq:Bayes} is based on prior population analyses leading to a statistical model
\begin{equation}
\Theta_i \sim p\big(\cdot \big|\hatthetaTV,\hatOmega\big) \label{eq:IIVModel}%
\end{equation}
with $\widehat{\theta^\text{TV}}$ denoting the estimates of the typical values (TV) and $\hat{\Omega}$ the estimated magnitude of inter-individual variability (IIV).
We assume in the following a normal distribution for $p(\cdot|\hatthetaTV,\hatOmega)$, which can typically be established via transformation, e.g., log-transformation in case of the log-normal distribution.
Note that it is possible to consider a covariate dependent model, but we restrict ourselves to a more fixed parametrization of the model parameters $\hatthetaTV$. The likelihood is based on structural and observational models
\begin{align}
\frac{\rd x_i}{\rd t}(t) &= f(x_i(t);\theta_i,d_i), \qquad x_i(0) =x_0(\theta_i)  \label{eq:ode} \\
h_i(t) &= h(x_i(t),\theta_i) \label{eq:observable} %
\end{align}
with state vector $x_i=x_i(t)$ (e.g., drug and neutrophil concentrations) and rates of change $f(x_i;\theta_i, d_i)$ for given doses $d_i$. 
The initial conditions $x_0(\theta_i)$ are given by the pre-treatment levels (e.g., baseline neutrophil concentration). 
A statistical model links the observables $h_{ij}(\theta_i) = h_i(t_{ij})$, i.e., model quantities that are measurable at time points $t_{ij}$, to observations $(t_{ij}, y_{ij})$ for $j=1,\ldots,n_i$ taking into account measurement errors and potential model misspecifications. This specifies the likelihood in \Eqref{eq:Bayes};  a common likelihood is defined by
\begin{equation}\label{eq:ResidualModel}
p\big(y_{ij} \big| \theta_i) = h_{ij}(\theta_i) + \epsilon_{ij};\qquad \epsilon_{ij} \sim_\iid \mathcal{N}(0,\Sigma)
\end{equation}
for $j=1,\dots,n_i$. For ease of notation, we do not explicitly state the dependence on $h_{ij}$ and $\Sigma$ in $p(y_{i}|\theta_{i})$.

In the simulation study, we considered the parameters in log-space and consider neutrophil-guided dosing, i.e., based on the observed neutrophil concentration. 
From the lowest neutrophil concentration within a treatment cycle (nadir) the neutropenia grade of a chemotherapy cycle is inferred, ranging from no neutropenia (grade 0), mild (grade 1), moderate (grade 2), to severe (grade 3) and life-threatening (grade 4).
We used an MIPD approach based on data assimilation (DA) to forecast the neutropenia time course, called DA-guided dosing, as presented in detail in \cite{Maier2020,Maier2020DARL}.
In short, a particle filter is used to approximate the posterior distribution \Eqref{eq:Bayes} at dose selection time points $t_c$ (start of cycle $c$), integrating data from the $i$-th patient up to $t_c$, i.e., $(t_{ij},y_{ij})$ with $t_{ij}\leq t_c$, via a sample approximation 
\begin{equation}
    p(\theta_{i}|y_{i1:t_c}) \approx  \sum_{m=1}^M w_{i}^{(m)} \mathbf{1}_{\{\theta_{i}^{(m)} = \;\theta_i\}}.
\end{equation}
It is based on an ensemble of weighted particles (samples) $\mathcal{E}_i=\lbrace \theta_{i}^{(m)},w_{i}^{(m)} : m=1,\ldots, M\rbrace$ comprising parameter values $\theta_{i}^{(m)}$ and importance weights $w_{i}^{(m)}$. 
The used particle filter included a resampling and rejuvenation step,
see \cite{Maier2020}.
Solving the structural and observational model \Eqref{eq:ode}+\Eqref{eq:observable} for each particle (sample) allows to compute the {\it a posteriori} probabilities of all neutropenia grades.
The optimal dose is determined by minimizing the weighted joint probability of life-threatening grade 4 ($P_{\text{grade}~4}$) and subtherapeutic grade 0 neutropenia ($P_{\text{grade}~0}$) for the next cycle, i.e.,
\begin{equation}\label{eq:OptGrade04}
    d^* = \underset{d \in \mathcal{D}}{\arg \min} \ \big\{ \lambda_1 P_{\text{grade}~4} + \lambda_2 P_{\text{grade}~0}\big\}.
\end{equation}
The weighting factors were chosen as $\lambda_1=2/3$ and $\lambda_2=1/3$ to penalize the risk to expose a patient to life-threatening infections (grade 4) more severely  than the risk to expose a patient to a subtherapeutic outcome (grade 0). The latter has been associated with reduced median overall survival \cite{DiMaio2005,DiMaio2006}.

\subsection{Simulation study framework}

\subsubsection*{Paclitaxel-induced neutropenia models}
\begin{figure}[H]
    \centering
    \includegraphics[width=\textwidth]{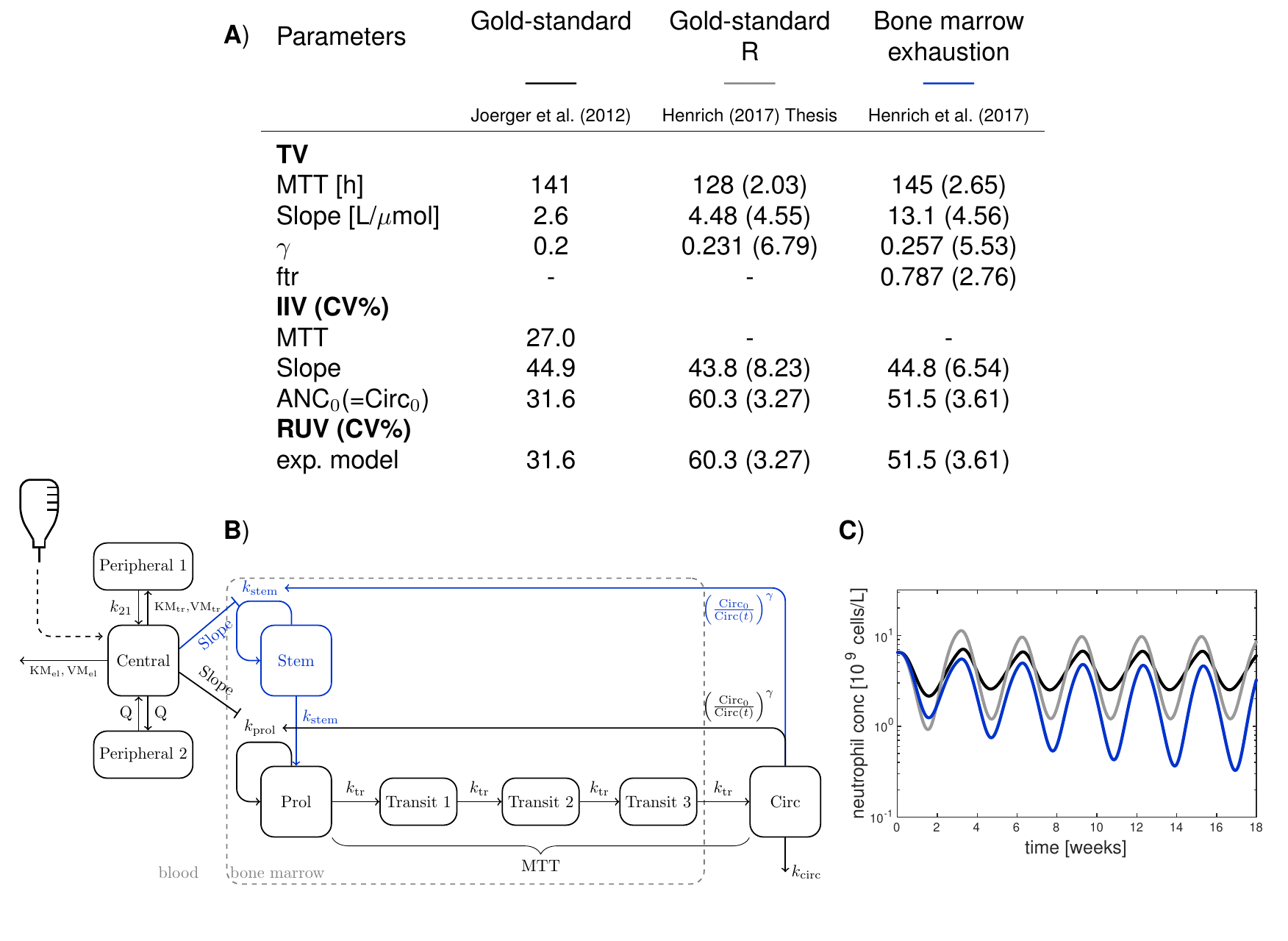}
    \caption{Paclitaxel-induced neutropenia models. \textbf{A)} Parameter estimates of the three considered models with relative standard errors in brackets (if available in literature). \textbf{B)} Schematic model representation with the PK model for paclitaxel (on the left); gold-standard PK/PD model in black and extension for cumulative neutropenia in black/blue. \textbf{C)} Time course of neutrophil concentrations for the typical patient in the CEPAC-TDM study  during six cycles of chemotherapy (three weeks each with drug administration on day 1) for the two model structures (black: gold/standard model trajectory; blue: extended gold-standard) and different parameter values (dark/light black)
\textbf{Abbreviations:} TV: typical values; MTT: mean transit time; Slope: linear drug effect parameter; $\gamma$: feedback exponent; ftr: fraction of input in the compartment of proliferating cells (Prol) via replication; IIV: inter-individual variability; ANC$_0$ = Circ$_0$: absolute neutrophil concentration at baseline or circulating neutrophils (Circ) at baseline; RUV: residual variability; BSA: body surface area; BILI: bilirubin concentration. More detailed information on the models is provided in Appendix \ref{sec:PKPDmodel}.}
    \label{fig:NeutropeniaModels}
\end{figure}

We considered a model for paclitaxel-induced neutropenia which, in ref.~\cite{Henrich2017}, had been investigated for model applicability, re-estimated and structurally modified after new patient data were observed in the clinical trial ``CESAR (Central European Society for Anticancer Research) study of Paclitaxel Therapeutic Drug Monitoring (CEPAC-TDM)'' [ClinicalTrials.gov Identifier: NCT01326767] \cite{Joerger2016}. 
A schematic representation of the models, corresponding parameter estimates and typical model predictions, is shown in Figure~\ref{fig:NeutropeniaModels}. 
In \ref{tab:PaclitaxelModelsALL} we list additional models proposed for paclitaxel-induced neutropenia, which illustrates the challenge of choosing a suitable model for MIPD in practice.
The initial model (hereafter \textit{gold-standard}) builds on the structure of the gold-standard model for chemotherapy-induced neutropenia \cite{Friberg2002} with parameter values estimated using a pooled data set of two prior studies \cite{Joerger2006,Joerger2007} including patients with ovarian cancer, non-small cell lung cancer (NSCLC) and various solid tumours \cite{Joerger2012}.
Paclitaxel was given either as monotherapy or in combination with carboplatin.
The CEPAC-TDM study included only NSCLC patients and paclitaxel was given in combination with carboplatin or cisplatin over six treatment cycles.
It was observed that the gold-standard model \cite{Joerger2012} overestimated the neutrophil concentration at later cycles since the model does not account for cumulative neutropenia, i.e., an aggravation of neutropenia over multiple cycles \cite{Henrich2017}, see Figure~\ref{fig:NeutropeniaModels}~C; a phenomenon that has been reported previously \cite{Huizing1997}.
The parameters were re-estimated (\textit{gold-standard~R}) based on the CEPAC-TDM data, and finally the structure was modified to account for bone marrow exhaustion (\textit{BME}), see Figure~\ref{fig:NeutropeniaModels}. 
Here, we focus our analyses on the more challenging PD models, while we considered the PK model to be given, with parameter values inferred previously based on the CEPAC-TDM study data \cite{Henrich2017}, see \ref{sec:PKPDmodel}.

\subsubsection*{Model adaptation scenarios} \label{sec:modelbias}

A model adaptation may be needed on the level of the 
structural model \eqref{eq:ode}+\eqref{eq:observable}, prior parameter distribution \eqref{eq:IIVModel} and/or likelihood \eqref{eq:ResidualModel}. 
In this article, we consider two types of scenarios where model adaptations may be beneficial:

\textit{Structural differences.} A divergence in the structural model, e.g., due to the manifestation of phenomena in the target patient population that have not been observed in the prior studies, is considered.
To study such a scenario, we used the BME model \cite{Henrich2017} to generate TDM data, while we used the gold-standard model \cite{Joerger2012} for MIPD. The latter model lacks the structural feature of cumulative neutropenia. 

\textit{Differences in parameters.} Differences in the parameter distribution e.g., the distributional assumption (normal, log-normal, etc.) as well as the estimated hyper-parameter values for a given distribution are potential examples.
Here, we only focus on the latter, i.e., the type of distribution is the same, but the hyper-parameters differ. 
To study parameter changes, we used the gold-standard~R model \cite{Henrich2017Thesis} 
to generate TDM data, while we used the gold-standard model \cite{Joerger2012} in MIPD. 
Both rely on the same structural model, while the parameter values of the former were re-estimated to the CEPAC-TDM data. 

We compared the performance of MIPD in the presence of structural or parameter bias to (i) MIPD based on an unbiased model (\textit{unbiased model} scenario), and (ii) standard dosing, i.e., 200\,mg/m$^2$ body surface area (BSA), including a dose reduction of 20\% if grade 4 was observed based on the neutrophil measurement at day 15 \cite{Henrich2017Thesis}.
%

\subsubsection*{TDM sampling scenarios}
The effect of a potential mismatch between prior model and the (new) data-generating process of MIPD depends on the amount of available TDM data per patient to adapt the model.
Therefore, we considered different TDM sampling schemes:
\begin{enumerate}
\item \textit{sparse sampling}: neutrophil measurements at day 1 and day 15 of each cycle (sampling design of CEPAC-TDM study).
\item \textit{intermediate sampling}: weekly neutrophil measurements (as in \cite{Joerger2007}).
\item \textit{rich sampling}: neutrophil measurements taken every third day.
\end{enumerate}
While the first two sampling schemes correspond to current clinical settings, the third mimics the prospective growing availability of point-of-care devices (e.g., HemoCue\textsuperscript{\textregistered} WBC Diff for measuring neutrophil counts \cite{Dunwoodie2018}), foreseeing richer sampling for monitoring patients.

\subsubsection*{Hierarchical Bayesian model}

\begin{figure}
    \centering
    \includegraphics[width=\textwidth]{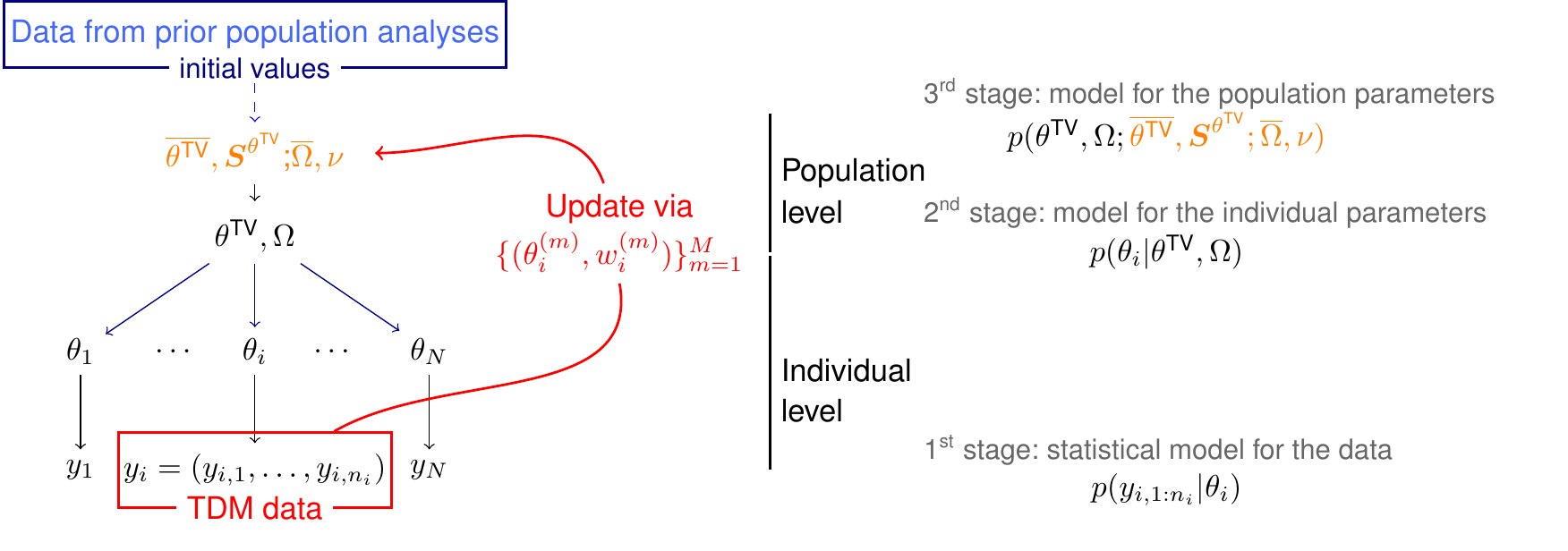}
    \caption{Hierarchical Bayesian model framework with separation between the inference on the individual level and the inference on the population level. The standard population analyses typically used to build the prior knowledge based on population estimates $\widehat{\theta^\text{TV}}$, $\hat{\Omega}$ is shown in dark blue. To update the prior population estimates, the population parameters are seen as random variables with parametric probability distributions parametrized with hyperparameters $\overline{\theta^\text{TV}}$, $\boldsymbol{S}^{\theta^\text{TV}}$; $\overline{\Omega}$, $\nu$. A sample representation of the individual posterior, $\lbrace\theta^{(m)}_{N+1}\rbrace_{m=1}^M$ is used to update the hyperparameters of the population parameter distributions (red arrow). On the right the corresponding probability distributions are given for the different levels of the hierarchical model.}
    \label{fig:HierarchicalBayes}
\end{figure}

To continuously update and learn population parameters, we considered additional hyper priors on the population parameters of the NLME models in \ref{sec:MIPDframework}.
The hierarchical structure of fully Bayesian population models thus comprises three stages \cite{Duffull2007,Wakefield1996}, see \Figref{fig:HierarchicalBayes}: 
(i) the statistical model for the TDM data given by \Eqref{eq:ResidualModel} describes the deviations between the individual model predictions and the observational data; 
(ii) the distributional assumption for inter-individual variability, \Eqref{eq:IIVModel} describes the differences between individuals; 
(iii) the distributional assumption for (hyper) population parameters, $p(\thetaTV), \  p(\Omega)$ describes the uncertainty in the population parameters.

Population analyses are typically performed in a frequentist NLME setting, reporting maximum likelihood estimates (MLE) of the population parameters jointly with their relative standard errors (RSE) or coefficients of variation (CV). 
This leaves the problem of how to determine suitable hyper prior distributions of the population parameters.
We considered a normal distribution for the typical parameters (on log-scale) and an inverse-Wishart for the variances, as suggested in  \cite{Gisleskog2002}.
We chose $\thetaTV$  to be normally distributed with mean $\overline{\thetaTV}$ identical to the MLE  $\hatthetaTV$ and variance $\mathbf{S}^{\thetaTV}$ identical to the (appropriately transformed) squared standard error $(\text{SE}_{\thetaTV})^2$, i.e. $\thetaTV\sim\mathcal{N}(\hatthetaTV,(\text{SE}_{\thetaTV})^2)$, see also \cite{Gelman2014}.
%
The inter-individual variability matrix $\Omega$ was assumed to be inverse-Wishart distributed $\mathcal{IW}(\Psi,\nu)$ and diagonal with parameters $\Psi,\nu$ such that the population estimate equaled the mean $\Psi/(\nu-n_{\Omega}-1)$, i.e., $\Psi=(\nu-n_{\Omega}-1)\hat{\Omega}$, where $n_{\Omega}$ is the number of random effect parameters and with degrees of freedom $\nu$ still to be chosen.
The distributions of the typical and variability values were assumed to be independent.

\subsection{Continued learning across patients}
To learn and improve a model across patients, the information provided by the patient-specific TDM data has to be included into the hierarchical model. 
In mathematical terms, we are interested in the marginal posterior:
\begin{equation}\label{eq:marginalposterior}
    p(\thetaTV,\Omega|y_{1:i}) \propto  \int p(\theta_{i},\thetaTV,\Omega|y_{1:i}) d\theta_{i}
\end{equation}
with the joint posterior
\begin{equation}\label{eq:jointposterior}
    p(\theta_{i},\thetaTV,\Omega|y_{1:i}) \propto p(y_{i}|\theta_{i})p(\theta_{i}|\thetaTV,\Omega)p(\thetaTV,\Omega|y_{1:i-1})
\end{equation}
determined from a full hierarchical Bayesian procedure.  
A sample approximation to the joint posterior \Eqref{eq:jointposterior} allows for a straightforward approximation of the marginal in \Eqref{eq:marginalposterior}.
In the context of particle filter-based inference, this can be realized by augmenting the particle state as well as the parameter space by the population parameters $(\theta^\text{TV},\Omega)$. 
This approach, however, has two major drawbacks: (i) it is computationally expensive and thus limits real-time inference during the patient's therapy; and (ii) direct access to the individual patient data is required to update the population parameters yet data protection laws and logistical reasons often prohibit this. These are major limitations for a practical application, in particular across different clinics.
\\

Therefore, we propose a two-level sequential Bayesian approach; it is based on previous ideas on Bayesian inference for meta-analyses \cite{Lunn2013,Hooten2019}. 
Importantly, this approach does not change the inference on the individual level (see Algorithm~\ref{alg:2LevelCLearning} for pseudo-code):
\begin{enumerate}
\item \textit{Individual level:} Estimate individual parameters of the $i$th patient 
\begin{equation*}
p(\theta_i|y_{i},\thetaTV,\Omega) \propto p(y_{i}|\theta_i) p(\theta_i|\thetaTV,\Omega)\,,\end{equation*}
e.g., using a particle filter, sampling importance resampling (SIR) or a Markov chain Monte Carlo (MCMC) approach \cite{Maier2020}.  A particle filter (`DA' for \textit{data assimilation} in the pseudo-code) is employed in our analysis, as it was shown to be best suited for the underlying setting \cite{Maier2020}.
This gives rise to a sample representation $\lbrace \theta_i^{(m)}, w_i^{(m)} : m=1,\dots,M \rbrace$ of the posterior $p(\theta_i|y_{i},\thetaTV,\Omega)$, summarizing the information provided by the data of the $i$th patient. 
This step is identical to the inference step in MIPD without continuous learning, see \Eqref{eq:Bayes}.

\item \textit{Population level:} Update population parameters by sampling iteratively from the joint posterior $p(\theta_{i},\thetaTV,\Omega|y_{1:i})$ via a Metropolis-Hastings-within-Gibbs sampling scheme \cite{Gelman2014,Lunn2013}, i.e., sampling from the full conditionals (see \ref{sec:deriv-full-conditionals} for a detailed derivation):

\begin{align}\label{eq:condthetaTV}
    p(\thetaTV|\theta_{i},\Omega,y_{1:i}) &\propto p(\theta_i|\thetaTV,\Omega) p(\thetaTV|\Omega,y_{1:i-1})\\[1em]
    \label{eq:condOmega}
    p(\Omega|\theta_{i},\thetaTV,y_{1:i}) &\propto p(\theta_i|\thetaTV,\Omega) p(\Omega|\thetaTV,y_{1:i-1})\\[1em]
    \label{eq:condTheta}
    p(\theta_{i}|\thetaTV,\Omega,y_{1:i}) &\propto p(y_{i}|\theta_i) p(\theta_i|\thetaTV,\Omega)
\end{align}

Sampling from Eq.~\ref{eq:condTheta} is achieved via a Metropolis-Hastings step, using as proposals the posterior samples generated on the individual level $\lbrace \theta_i^{(m)} \rbrace_{m=1}^M$, which are drawn according to weights $w_i^{(m)}$.

To ensure that sampling in Eqs.~\ref{eq:condthetaTV}-\ref{eq:condOmega} can be performed in closed-form, at the end of  assimilating data of the $i$-th patient, a parametric approximation by a normal-inverse-Wishart distribution is used,
\begin{equation}
\label{eq:hyperprior-approx}
p(\thetaTV,\Omega|y_{1:i}) \approx \mathcal{N}(\overline{\thetaTV}_{i}, \mathbf{S}^{\thetaTV}_{i}) \otimes \mathcal{IW}((\nu_i-n_\Omega-1)\overline{\Omega}_{i},\nu_{i}),
\end{equation}

with hyperparameters as stated in Algorithm~\ref{alg:2LevelCLearning}. 
Then, sampling in Eqs.~\ref{eq:condthetaTV}-\ref{eq:condOmega} corresponds to sampling from a normal  and inverse-Wishart distribution, respectively, see \ref{sec:CLapproach}.

Importantly, through the parametric approximation in Eq.~\ref{eq:hyperprior-approx}, 
$y_{1:i-1}$ are represented implicitly via the updated priors for the $i$-th step, while $y_i$ enters implicitly through the sample approximation. In no case, the original patient data is needed. 


    \end{enumerate}

\begin{algorithm}
\caption{Two-level sequential hierarchical Bayesian learning in MIPD}\label{alg:2LevelCLearning}
\begin{algorithmic}[1]
\State Input: $\hatthetaTV, \text{SE}_{\hatthetaTV} , \hat{\Omega},\nu_0$, ($y_{i,1:n_i}$ only for individual level)
\State Set hyper prior parameters $\overline{\thetaTV_0}:= \hatthetaTV,\mathbf{S}_0^{\thetaTV}:=(\text{SE}_{\hatthetaTV})^2,\bar{\Omega}_0:=\hat{\Omega},\nu_0$
\For {$i=1:N_\text{TDM}$}
    \LineComment{Individual level}
    \State Initialize particle ensemble $\lbrace (\theta_{i0}^{(m)},x_{i0}^{(m)},w_{i0}^{(m)}) \rbrace_{m=1}^M$ based on $p(\theta|\overline{\thetaTV_{i-1}},\bar{\Omega}_{i-1})$
    \For {$j=1:n_i$}
        \State $\lbrace (\theta_{ij}^{(m)},x_{ij}^{(m)},w_{ij}^{(m)}) \rbrace_{m=1}^M \leftarrow \text{DA}\big(y_{ij},\lbrace (\theta_{ij-1}^{(m)},x_{ij-1}^{(m)},w_{ij-1}^{(m)}) \rbrace_{m=1}^M\big)$
    \EndFor
    \LineComment{Population level} 
    \State Initialize Markov chain ${\thetaTV_i}^{(0)}=\overline{\thetaTV_{i-1}},\Omega^{(0)}_i=\bar{\Omega}_{i-1}$ and $\theta_i^{(0)}$ sampled from $p(\theta|\overline{\thetaTV_{i-1}},\bar{\Omega}_{i-1})$
    \For {$l=1:L$}
    \LineCommentg{Gibbs sampling part}
    \State Draw $\thetaTVl_i$ from $p(\thetaTV|\theta_i^{(l-1)},\Omega_i^{(l-1)})$ \Comment{\Eqref{eq:condthetaTV}}
    \State Draw $\Omega^{(l)}_i$ from $p(\Omega|\theta_i^{(l-1)},\thetaTVl)$
    \Comment{\Eqref{eq:condOmega}}
    \LineCommentg{Metropolis-Hastings part}
    \State Draw proposal $\theta^{*(l)}_i$ from $\lbrace \theta_{in_i}^{(m)} \rbrace_{m=1}^M$ according to $\lbrace w_{in_i}^{(m)} \rbrace_{m=1}^M$ and add rejuvenation
    \State Accept proposal with probability
        \begin{equation}\label{eq:AccRatio}
        \alpha = \frac{p(\theta^{*(l)}_i|{\thetaTVl_i},\Omega_i^{(l)})/p(\theta_i^*|\overline{\thetaTV_{i-1}},\bar{\Omega}_{i-1})}{p(\theta_i^{(l-1)}|\thetaTVl_i,\Omega_i^{(l)})/p(\theta_i^{(l-1)}|\overline{\thetaTV_{i-1}},\bar{\Omega}_{i-1})}\,,
    \end{equation}
    %

    \EndFor
    \State Parametric approximations of posterior (hyper prior for next individual): 
    \State $p(\thetaTV|y_{1:i})\approx \mathcal{N}(\overline{\thetaTV_i},\mathbf{S}^{\thetaTV}_i)$ with \begin{equation}\overline{\thetaTV_i} = \frac{1}{L}\sum_{l=1}^L \thetaTVl_i, \quad  \mathbf{S}^{\thetaTV}_i=\frac{1}{L-1}\sum_{l=1}^L (\thetaTVl_i-\overline{\thetaTV_i})(\thetaTVl_i-\overline{\thetaTV_i})^T\end{equation}
    \State $p(\Omega|y_{1:i}) \approx \mathcal{IW}((\nu_i-n_\Omega-1)\bar{\Omega}_i,\nu_i)$ with 
    \begin{equation}{\bar{\Omega}_i}=\frac{1}{L}\sum_{l=1}^L {\Omega_i^{(l)}},\, \quad \nu_i =\nu_{i-1}+1\end{equation}
\EndFor
\end{algorithmic}
\end{algorithm}

\noindent
The continued learning approach was sequentially applied to $N_\text{TDM}=100$ virtual patients with available TDM data over six treatment cycles depending on the considered sampling scheme.
For the analysis, the continuous learning approach was repeated $10$ times to account for statistical variability in the individual patient parameters considered for the update. 
To demonstrate the effect of a mismatch between model and data-generating process, we also applied MIPD alone without continued learning (DA-guided dosing).
On the individual level, model parameters $(\text{MTT}, \text{Slope}, \text{ANC}_0)^T$ were estimated. 
We restricted the population updates to `MTT' and `Slope' as for `ANC$_0$' the baseline method B2 described in \cite{Dansirikul2008} was used, i.e., no typical parameter was estimated but the baseline value was used to initialize the (empirical Bayes) prior.
In addition, we consider a setting which includes $\gamma$ in the individual level inference as the value differs across the models.
In this study we neither estimated the residual variability $\sigma$ on the individual level nor on the population level.
The values for $\sigma$ used to generate the TDM data, however, differed from those that the models in the `imperfect model scenarios' assume.
The considered hyper priors, i.e., the distributional assumptions for the population parameters, are summarized in \Tabref{tab:Hyperpriors}.
Since no relative standard errors are available for the gold-standard model \cite{Joerger2012}, the values reported in \cite{Joerger2007} (one of the two pooled studies) were chosen as conservative choice. 
The degrees of freedom $\nu$ was chosen here to balance confidence in the estimated value while still enabling adaptation.
The simulation study was performed in MATLAB 2019b and the code is available under \url{https://doi.org/10.5281/zenodo.39670}.

\begin{table}[H]
\caption{Hyper priors for the gold-standard model used in the simulation study.}
\label{tab:Hyperpriors}
\begin{center}
\begin{tabular}{lcc}
Parameter & distribution & hyperparameters \\
\hline
\textbf{TV parameters} 	&		&		       \\
log(MTT) 				& $\mathcal{N}$	& $\overline{\thetaTV}_0=\log(2.6)$ , $\mathbf{S}_0^{\thetaTV}= 0.0013$	 	\\
log(Slope)				& $\mathcal{N}$	& $\overline{\thetaTV}_0=\log(141)$	, $\mathbf{S}_0^{\thetaTV}=0.016$	\\
	
\textbf{IIV parameters } 	& 		&			\\	
$\omega^2_\text{MTT}$	 				& $\mathcal{IW}$	& $\Psi_\text{MTT} = 0.6561=(12-2-1)0.0729$, $\nu_0= 12$	\\
$\omega^2_\text{Slope}$					& $\mathcal{IW}$	& $\Psi_\text{Slope} = 1.8144=(12-2-1)0.2016$,	$\nu_0=12$	\\
\end{tabular}
\end{center}
\end{table}

\section{RESULTS}\label{sec:Results}

\subsection{Current MIPD approaches may not be beneficial in the presence of model bias}

\begin{figure}[H]
\includegraphics[width=\textwidth]{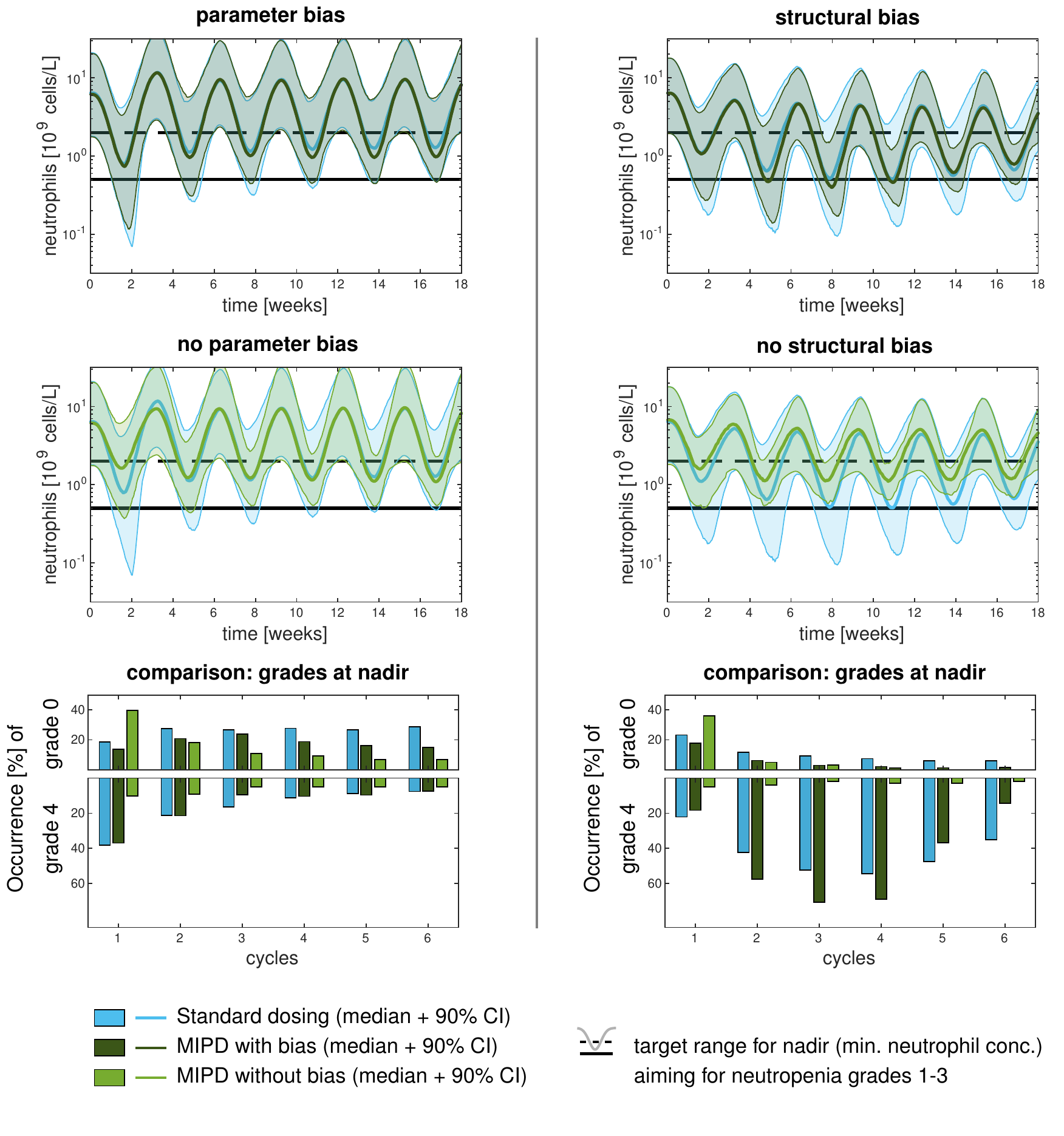}
\caption{\textbf{MIPD in presence of a mismatch between model and data-generating process}. Neutropenia time-courses were simulated for $N_\text{TDM}=1000$ virtual patients using MIPD approaches with different underlying models and the standard dosing approach for paclitaxel ($200\,\text{mg}/\text{m}^2$ with 20\% reduction if grade 4 neutropenia is observed in the previous cycle). 
TDM data were simulated for day 1 \& day 15 of each cycle (sparse sampling scenario), including residual variability. 
\textbf{Left column}: TDM data were generated using the gold-standard~R model. \textbf{Right column}: TDM data was generated using the bone marrow exhaustion (BME) model \cite{Henrich2017}.
\textbf{Top panel}: `imperfect model scenario', i.e., the MIPD approach uses the gold-standard model, while the data was generated with the gold-standard~R model (left) or BME model (right). \textbf{Middle panel}: `perfect model scenario', i.e., the MIPD approach uses also the data-generating model. The median time-course is shown along with its 90\% confidence interval (CI).
\textbf{Bottom panel}: occurrence of life-threatening grade 4 neutropenia and subtherapeutic grade 0 neutropenia for the different scenarious.
Note that grade 4 neutropenia is penalized more ($\lambda_1=2/3$) compared to grade 0 neutropenia ($\lambda_2=1/3$) in \Eqref{eq:OptGrade04}. This is accounted for in the scale of the bottom panels, which allows to interpret the length of the total bars.
}
\label{fig:MIPDwithwoModelbias}
\end{figure}

For a performance analysis, we generated TDM data (including residual variability) on day 1 \& day 15 of each cycle (sparse sampling as in the CEPAC-TDM study). 
\Figref{fig:MIPDwithwoModelbias} illustrates the performance of MIPD with/without model bias in comparison to standard dosing (median and 90\% confidence intervals (CIs)).

The left column illustrates the scenario of parameter deviation, i.e., the structural model and the class of prior distributions is identical to the data-generating process, but hyper-parameter values differ. 
In this case, MIPD performed comparably to standard dosing (top left, median trajectory and 90\% CI), also in terms of occurrence of grade 4 \& 0 neutropenia (bottom left). 
For reference, in the corresponding model scenario without mismatch, the MIPD approach clearly reduced the occurrence of grade 4 \& 0 (bottom and middle panel). 
It is worth mentioning that the CIs in all panels showed a certain `skewness' towards higher neutrophil concentrations (lower neutropenia grade), since grade 4 is penalized more strongly than grade 0 in \Eqref{eq:OptGrade04} .

The right column illustrates the more challenging scenario of structural changes, i.e., a model structure differing from the data-generating process.
Of note, in this case both standard dosing and MIPD performed much worse than in the parameter bias scenario (bottom panels). 
In 3 out of 6 cycles (cycles 2-4), MIPD resulted in even larger occurrences of grade 4 compared to standard dosing. The gold-standard model underestimated the drug effect on neutrophil concentrations (see \Figref{fig:NeutropeniaModels}) and too high doses were selected, especially in presence of cumulative neutropenia. 
Despite relying on an inappropriate structural model, DA was able to correct this initial mismatch on the parameter level over the course of a patient's therapy by integrating TDM data, which lead to a decrease in incidence of grade 4 neutropenia in later cycles. 
For reference, in the corresponding unbiased model scenario, the MIPD approach clearly and very quickly reduced the occurrence of grade 4 \& 0 (bottom and middle panel). 
In comparison to the parameter adaptation scenario, the occurrence of grade 4 \& 0 neutropenia was even further decreased, which might be related to the smaller RUV parameter, see \Figref{fig:NeutropeniaModels}~A.

In summary, if the underlying model is not consistent with the observational data, MIPD might not be beneficial compared to standard dosing that solely relies on TDM data (`model-free'). 
As outlined in the introduction, necessary model adaptations can be expected if MIPD is applied in clinical routine, and therefore, the top panels might better reflect clinical reality than the middle panels. 
Here, model adaptation during clinical practice is necessary. Most MIPD approaches, however, do not exploit the wealth of TDM data used during MIPD to learn and update their models. 

\subsection{Continued learning MIPD can adapt parameters---depending on the sampling scheme}

\begin{figure}[H]
    \centering
    \includegraphics[width=\linewidth]{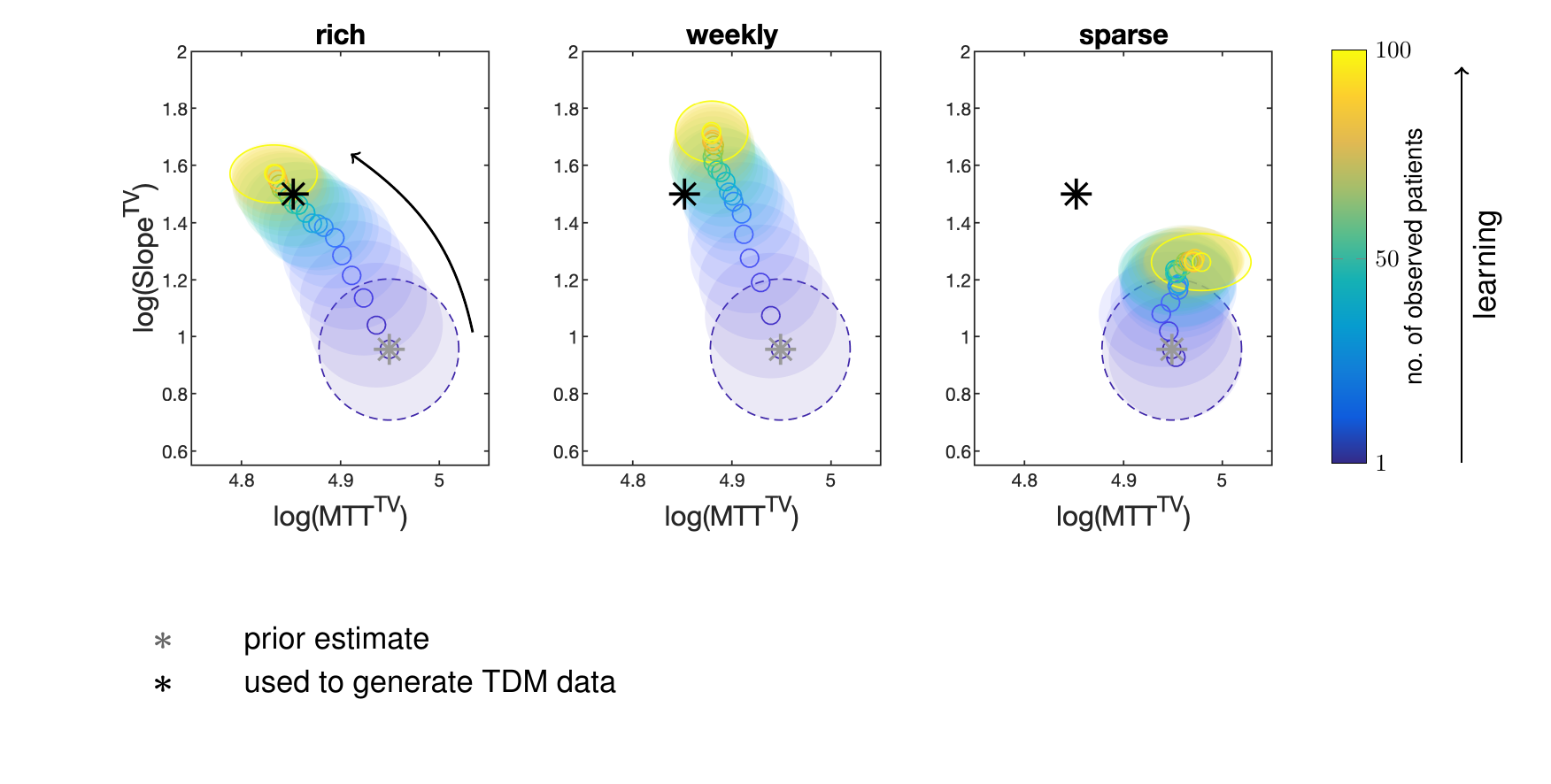}
    \caption{Comparison of the sequential updates of the hyper prior for the typical MTT and Slope value for different TDM scenarios. Grey star: prior estimate of the hyper-parameters; black star: true hyper-parameters, i.e., the values used to generate the TDM data. Sparse sampling consisted of measurements on day 1 \& 15 of each cycle, weekly sampling corresponds to an intermediate data situation and for rich sampling it is assumed that neutrophils are monitored every third day. Mean (circle) and 95\% highest posterior density (shaded ellipse) are shown.}
    \label{fig:NAP_logthetaTV1_20}
\end{figure}

The proposed continued learning framework was able to adapt the prior parameter distribution across TDM patients. 
Figure~\ref{fig:NAP_logthetaTV1_20} illustrates the sequential updates of the proposed framework for the posterior distributions of the typical parameters of `Slope' and `MTT' across $100$ patients for different sampling schemes. 
For the rich sampling scenario (left), the posterior---95\% highest posterior density (HPD) area---evolved over the number of observed patients, moving away from the prior estimate (grey star) towards the value used to generate the data (black star).
As more patients were observed, uncertainty about the typical `Slope' and `MTT' parameters decreased, as indicated by the decreasing size of the HPD area.
Thus, the proposed framework successfully learned the typical values underlying the TDM data from sample representations of the posterior on the individual level.
Note that the parameters $\gamma$ and $\sigma$ were not estimated although different values were used to generate the data, which has the effect of introducing an additional bias.
The results including $\gamma$ on the individual level inference are shown in \ref{fig:NAP_mvlogthetaTV1_100_gamma}.

The extent to which the continuous learning framework could counteract a parameter mismatch depended on the sampling scheme, see \Figref{fig:NAP_logthetaTV1_20} (middle and right panel). 
For the intermediate scheme, the posterior distribution moved towards the parameter values used to generate the data. 
A final parameter bias, however, remained. A potential reason could be parameter identifiability. 
To assess practical identifiability, we investigated the log-likelihood and log-posterior on the individual patient levels, see \ref{fig:Identifiability}. 
To exclude the possibility of unfavorably chosen sampling time points in the intermediate scheme (weekly), we performed an optimal design analysis, see \ref{sec:SupplementAnalyses}.
For the sparse sampling scheme, the TDM data were not sufficient to adapt the model appropriately. 
Yet in the context of the rich sampling the data was indeed informative enough to move away from the (biased) prior estimate towards the data-generating value, resolving the practical unidentifiability.

\subsection{Continued learning in MIPD can substantially improve therapy outcome even for structural changes}

\begin{figure}[H]
\includegraphics[scale=1]{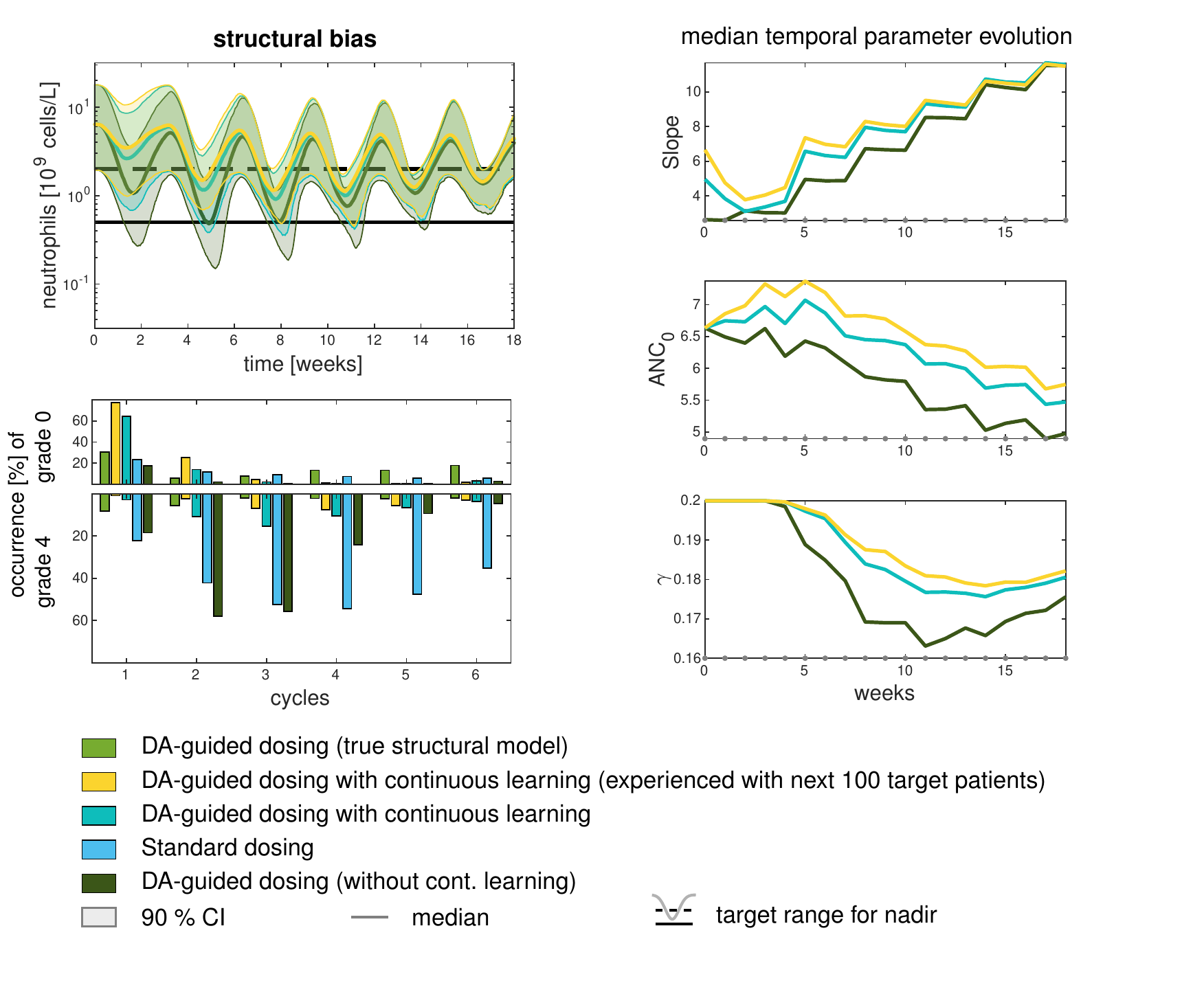}
\caption{Sequential DA allows for temporal parameter changes within the course of a patient's therapy. The TDM data were generated using the BME model and the results are shown for the intermediate sampling scheme (sampling timepoints are indicated as grey dots in the right panels). The median temporal parameter evolution over the course of therapy was computed across all virtual patients ($N_\text{TDM}=10\cdot 100$)}
\label{fig:TemporalChanges}
\end{figure}

Finally, we investigated the effects of continued learning of population parameters on MIPD. 
Here, we show only the more challenging structural bias scenario (a different model used for data generation vs.\ inference in MIPD); 
for the parameter bias scenario, see \ref{fig:ImpactonMIPD}.
The performance of the different approaches is compared in \Figref{fig:TemporalChanges}. Standard dosing and DA-guided dosing are set-up analogously in \Figref{fig:MIPDwithwoModelbias}, except for the fact they are trained with the intermediate sampling design. It becomes clear is that continued learning is significantly more effective with more TDM data (see above).
We also considered uncertainty with respect to the parameter $\gamma$.

DA-guided dosing was also able to adjust to some extent to cumulative neutropenia over time, see \Figref{fig:TemporalChanges} (dark green).
The `Slope' parameter increased, while parameters `Circ$_0$' and $\gamma$ decreased over the course of the individual therapy, leading to a decrease in occurrence of grade 4 after cycle 3 and a substantial decrease in outcome variability.
Effectively, when considering the data points one at a time, the sequential DA framework allowed to account for changes in the parameters over time---a potentially very beneficial property (e.g., in disease progression). While this might be very desirable for MIPD at the individual patient level, it could be misleading when learning across patients. When the final parameter estimate (after six cycles) was used to update the population parameter (Slope$^\text{TV}$, MTT$^\text{TV}$), this introduced a bias for the first cycle of the next patient, resulting in high occurrence of grade 0 for the first cycle (\Figref{fig:TemporalChanges} bottom left). 
Continued learning was considered across the first 100 patients (blue-green) as well as for some second 100 patients (yellow) after learning from the first $100$ patients. 
It can be observed that the typical `Slope' parameter increased (green vs. blue-green vs. yellow) as it was continuously learned across patients (initial Slope value at $t=0$ in the top right panel).

A major improvement was observed for the continued learning MIPD approach, which reduced the occurrence of grade 4 substantially across all cycles compared to DA-guided dosing or standard dosing.
The results for the rich sampling scenario were comparable (\ref{fig:MB3gammaRICH}); for the sparse sampling scheme, however, the benefits were not so clear (\ref{fig:MB3gammaSPARSE}).

\section{DISCUSSION}\label{sec:Discussion}
We proposed a sequential hierarchical Bayesian approach to update the population parameters using posterior samples as a means to exchange information with every treated patient so that the model better reflects the target patient population.
We showed that the approach allowed to successfully learn the underlying population parameters of the PK/PD model used to generate the patient data. 
It is important to note, however, that the results depend on the sampling scheme.
In addition, we showed that continued learning has potential to improve MIPD even in presence of structural changes, again depending on the informativeness of available TDM data.

The proposed approach has two levels and allows to learn sequentially over patients without using patient data on the population level.
Thus, the patient data itself does not need to be stored or shared across centers, which is a crucial advantage compared to pooling approaches \cite{Hughes2020}.
Thus, the approach builds a basis to develop more informed models, integrating an ever-growing database potentially better reflecting rare covariates.
%
The initial model used to start the continued learning approach could be selected using a retrospective external evaluation based on historical data from the intended clinical setting \cite{Zhao2013}.
Model selection/model averaging approaches do not adapt/improve the underlying model across patients; 
the {\it a priori} forecast remains the same for all patients (based on the covariates).
In addition, in their general form these approaches are implemented in conjunction with MAP estimation which provides potentially biased predictions in context of nonlinear models \cite{Maier2020}.
The proposed DA-guided dosing also naturally extents to model averaging and this extension has been considered (on the individual inference level) previously in the context of Bayesian therapy forecasting \cite{Albers2017}.\\

In the context of cytotoxic chemotherapy with neutropenia as dose-limiting toxicity, we show how population-based PK/PD models can be transferred to a different clinical target population, which is often a crucial application hurdle of MIPD in clinical routine.
We showed that model misspecification might severely impact MIPD, and therefore, models should be adapted to the target patient population.
The used DA-guided dosing proved to be able to adapt the model to some extent, but only improved MIPD at later cycles, when a certain amount of TDM data was collected. 
This is a consequence of assimilating data sequentially, thereby allowing to account for temporal changes in the parameters and thus adapting the gold-standard model to some extent to cumulative neutropenia.

The presented analysis revealed that an important aspect for practical implementation is to critically assess the quality of the inference on the individual level.
The dependence on the sampling design clearly showed that more research is necessary and that caution is needed when updating models based on real-world data, as was demonstrated for time-dependent parameters.
With the prospect of novel digital health care devices, e.g. point-of-care devices, more frequent monitoring could become clinical reality.
Real-world data is currently underutilized \cite{Tyson2020} but has great potential to improve MIPD as shown in this study.
The current approach is limited to misspecifications or population shifts on the structural model parameter level.
A model change on the structural level, e.g., accounting for cumulative neutropenia, could be corrected for to some extent on the level of parameters.
An important extension in the future would be to also estimate the RUV parameter $\sigma$, as an increased error in measurement precision or reporting can be expected in clinical routine compared to e.g.,~controlled clinical study settings.
Currently, the IIV parameters $\omega^2$ captured the increased RUV of the data to some extent, which, however, also increased the uncertainty on the individual level, see \ref{fig:mvGMH_Omega}.\\

The approach of a learning model (as coined by \cite{Keizer2018}) for MIPD could be beneficial not only in clinical practice, but also during drug development, where new (clinical) study data are generated continuously and should be integrated into previously developed models \cite{Lalonde2007}.
The study is an important step towards building the underlying models of MIPD on a growing database and thus make MIPD fit-for-purpose in everyday therapeutic use.

\section*{Acknowledgement}

C.M. kindly acknowledges financial support from the Graduate Research Training Program PharMetrX: Pharmacometrics \& Computational Disease Modelling, Berlin/Potsdam, Germany. This research has been partially funded by 
Deutsche Forschungsgemeinschaft (DFG) - SFB1294/1 - 318763901.
Fruitful discussions with Sven Mensing (AbbVie, Germany), Alexandra Carpentier (Otto-von-Guericke-Universitaet Magdeburg), Sebastian Reich (University of Potsdam, University of Reading) and David Albers (University of Colorado) are kindly acknowledged.

\subsubsection*{Conflict of Interest/Disclosure}
CK and WH report research grants from an industry consortium (AbbVie Deutschland GmbH \& Co. KG, AstraZeneca, Boehringer Ingelheim Pharma GmbH \& Co. KG, Gr\"unenthal GmbH, F. Hoffmann-La Roche Ltd., Merck KGaA and Sanofi) for the PharMetrX PhD program. In addition CK reports research grants from the Innovative Medicines Initiative-Joint Undertaking (`DDMoRe'), the European Commission within the Horizon 2020 framework programme (`FAIR'), and Diurnal Ltd. All other authors declare no competing interests for this work.

\subsubsection*{Funding information}
\begin{itemize}
\item Graduate Research Training Program PharMetrX: Pharmacometrics \& Computational Disease Modelling, Berlin/Potsdam, Germany,
\item Deutsche Forschungsgemeinschaft (DFG) - SFB1294/1 - 318763901.
\end{itemize}

\section*{Author Contributions}

C.M., N.H., J.dW., C.K., W.H., designed research, C.M. mainly performed the research, C.M., N.H., J.dW., C.K., W.H.,  analysed data and wrote the manuscript.
\bibliographystyle{cptpspref}
\bibliography{Modelbias_article}

\appendix

\section{PK/PD model for paclitaxel-induced cumulative neutropenia}\label{sec:PKPDmodel}
We employed published models describing the pharmacokinetics (PK) of paclitaxel as well as the side effects of the drug on the hematopoietic system (pharmacodynamics, PD). 
Some of the models were previously described in the supplementary material in \cite{Maier2020}; here, we repeat this description for convenience.

\subsection{Paclitaxel PK model}
Paclitaxel is a widely used anticancer drug, approved for the treatment of advanced non-small cell lung cancer \cite{Belani2005}. 
Paclitaxel PK was previously described by a three-compartment model with nonlinear distribution to the first peripheral compartment and nonlinear elimination \cite{Joerger2012}.
For all our analyses, we used the (re-estimated) parameter values by \cite{Henrich2017}, see \ref{tab:PaclitaxelParameterEstimates}. 
The PK model includes the covariate model
\begin{align*}
\mathrm{VM}_\mathrm{EL,TV,i} = \ &\mathrm{VM}_\mathrm{EL,pop} \cdot \Big(\frac{\mathrm{BSA}_i}{1.8 m^2}\Big)^{\theta_{\mathrm{VM}_\mathrm{EL}\text{-}\mathrm{BSA}}} \cdot \Big(\theta_{\mathrm{VM}_\mathrm{EL}\text{-}\mathrm{SEX}}\Big)^\mathrm{SEX_i}\\
& \cdot \Big(\frac{\mathrm{AGE}_i}{56 \mathrm{years}}\Big)^{\theta_{\mathrm{VM}_\mathrm{EL}\text{-}\mathrm{AGE}}} \cdot \Big(\frac{\mathrm{BILI}_i}{7 \mu \mathrm{mol}/\L}\Big)^{\theta_{\mathrm{VM}_\mathrm{EL}\text{-}\mathrm{BILI}}}\,,
\end{align*}
where $\mathrm{BSA}$ denotes the body surface area, $\mathrm{SEX}$ the patient's gender ($0/1$ for female/male),  $\mathrm{AGE}$ the patient's age and $\mathrm{BILI}$ the bilirubin concentration.
In addition to inter-individual variability and residual variability, inter-occasion variability was included on the two parameters $V_1$ and $\mathrm{VM}_\mathrm{EL}$, with an occasion being defined as one chemotherapeutic cycle $c$,
\begin{equation*}
\theta_{i,c} = \theta_\mathrm{TV,i} \cdot e^{\eta_i + \kappa_{i,c}}\,, \qquad \eta_i \sim_\text{iid} \mathcal{N}(0,\Omega), \ \kappa_{i,c} \sim_\text{iid} \mathcal{N}(0,\Pi)\,.
\end{equation*}

Paclitaxel PK was described by the following system of ordinary differential equations (ODEs):

\begin{align*}
\frac{d \Cent}{dt} &= u(t) - \frac{C_1\cdot \mathrm{VM}_\mathrm{EL}}{\mathrm{KM}_\mathrm{EL}+C_1} + k_{21} \text{Per1} - \frac{C_1\cdot \mathrm{VM}_\mathrm{TR}}{\mathrm{KM}_\mathrm{TR}+C_1}  + k_{31}\text{Per2} - k_{13}\Cent\,, &\Cent(0) = 0 \\
\frac{d \text{Per1}}{dt} &= \frac{C_1\cdot \mathrm{VM}_\mathrm{TR}}{\mathrm{KM}_\mathrm{TR}+C_1} - k_{21}\text{Per1}\,, &\text{Per1}(0)=0\\
\frac{d \text{Per2}}{dt} &= k_{13}\Cent - k_{31}\text{Per2}\,, &\text{Per2}(0)=0\\
\end{align*}
where $C_1(t) = \Cent/V_1$.

\begin{table}
\begin{center}
\begin{tabular}{lll}
\begin{tabular}{lll}
\hline
\multicolumn{3}{c}{structural submodel}\\
\hline
$\mathrm{V_1}$ 				& 	10.8 		& $[\L]$\\
$\mathrm{V_3}$ 				& 	301 		& $[\L]$\\
$\mathrm{KM}_\mathrm{EL}$ 		& 	0.667 	& $[\mu \M]$\\
$\mathrm{VM}_\mathrm{EL,pop}$ 	& 	35.9  	& $[\mu \mathrm{mol}/h]$\\
$\mathrm{KM}_\mathrm{TR}$ 		& 	1.44	 	& $[\mu \M]$\\
$\mathrm{VM}_\mathrm{TR}$ 		& 	175   	& $[\mu \mathrm{mol}/h]$\\
$k_{21}$						&	1.12 	 	& $[1/h]$\\
$Q$							&	16.8		& $[1/h]$ \\
\hline
\multicolumn{3}{c}{covariate submodel}\\
\hline
$\theta_{\mathrm{VM}_\mathrm{EL}\text{-}\mathrm{BSA}}$ &  1.14 &\\
$\theta_{\mathrm{VM}_\mathrm{EL}\text{-}\mathrm{SEX}}$ &  1.07 &\\
$\theta_{\mathrm{VM}_\mathrm{EL}\text{-}\mathrm{AGE}}$ &  -0.447 &\\
$\theta_{\mathrm{VM}_\mathrm{EL}\text{-}\mathrm{BILI}}$  &  -0.0942 &\\
\end{tabular}
\hphantom{xxxx}&
&
\begin{tabular}{lll}
\hline
\multicolumn{3}{c}{statistical submodel IIV}\\
\hline
\\
$\omega^2_{V_3}$ 					&  0.1639 &\\
\\
$\omega^2_{\mathrm{VM}_\mathrm{EL}}$ &  0.0253 &\\
$\omega^2_{\mathrm{KM}_\mathrm{TR}}$ &  0.3885 &\\
$\omega^2_{\mathrm{VM}_\mathrm{TR}}$ &  0.077 &\\
$\omega^2_{k_{21}}$ 					&   0.008  &\\
$\omega^2_{Q}$ 					&  0.1660 &\\

\hline
\multicolumn{3}{c}{statistical submodel IOV}\\
\hline
$\pi^2_{V_1}$ 					&   0.1391 &\\
$\pi^2_{\mathrm{VM}_\mathrm{EL}}$ 					&  0.0231&\\
\hline
\multicolumn{3}{c}{statistical submodel RV}\\
\hline
$\sigma^2$ 						& 0.0317& \\
\end{tabular}
\end{tabular}
\end{center}
\caption{Re-estimated PK parameter estimates of the previously published PK model \cite{Joerger2012} for the anticancer drug paclitaxel \cite{Henrich2017}.}
\label{tab:PaclitaxelParameterEstimates}
\end{table}

\subsection{Gold-standard model structure}
The gold-standard models in the main text are based on the structure of the chemotherapy-induced neutropenia model introduced in~\cite{Friberg2002}.
It describes the effect of various anticancer drugs on the hematopoietic system.
More specifically, a linear drug effect is assumed on the proliferation rate $k_\text{prol}$ of proliferating cells in the bone marrow (Prol).
Proliferating cells in the bone marrow differentiate over multiple progenitor cell stages to neutrophils, which are released into the systemic circulation.
Therefore, a decreased proliferation rate in the bone marrow becomes apparent with a certain delay in the amount of neutrophils.
Low neutrophil concentrations in turn lead to an increased proliferation rate, which restores normal neutrophil levels. 
A schematic representation of the model along with parameter estimates for paclitaxel is provided in Figure~1 in the main text.
In \ref{tab:PaclitaxelModelsALL} we provide multiple models which were all proposed for paclitaxel-induced neutropenia based on different or even on the same patient population.
This shows the difficulty associated with the selection of a model to use for MIPD.

\newpage

\begin{sidewaystable}
\begin{center}
\begin{adjustbox}{max width=\textwidth}
\begin{tabular}{lccccccc}
\multirow{2}{*}{Parameter} & Friberg et al. & Kloft et al. & Hansson et al. & Joerger et al. & Joerger et al. & Henrich et al & Henrich et al. \\
& (2002) & (2006) & (2010) & (2007) & (2012) & (2017) & (2017) \\
\hline
N (Ids)					&	45		&	45			& 	45		& 104 		&	104	& 366		& 	366	\\
n (samples)				&	530		&	530			&	523		&			&	314	& 3274		& 	3274	\\
cycles 					&	3		&	3 (1-18)		& 3 (1-11)		& (1-10)		&	1	& 6			& 	6	\\
structure 					&	G-S		&	G-S			& G-S		& G-S		&	G-S	& G-S		& 	BME	\\
\textbf{TV} 				&			&				&			&			&		&			& 		\\
ANC$_0$ [$10^9$ cells/L]		& 5.20 (3.6)	& 5.40 (7.2)		& 5.61 (9.4)	& 4.35 (1.5-9.4)* & -		& 6.48*		& 6.48*  	\\
MMT [h]					& 127 (2.1)	& 126 (4.2)		& 154 (4.4)	& 141 (3.7)	& 141	& 128 (2.03)	& 145  (2.65)	\\
Slope [L/$\mu$mol]			& 2.21 (4.5)	& 2.8 (13)			& 3.48 (8.1)	& 2.08$^\dagger$ (12.5)	& 2.6	 	& 4.48 (4.55)	& 13.1 (4.56)	\\
$\gamma$				& 0.230 (2.8) 	& 0.223			& 0.270 (5.9)	& 0.26 (7.5)	& 0.2		& 0.231 (6.79)	& 0.257 (5.53)	\\	
ftr						& -			& -				& -			& -			&	-	&	-		& 0.787 (2.76)	\\
\textbf{IIV
(CV\%)} 					& 			&				&			&			&		&			& \\	
MMT	 					& 18	(30)		& 17	 (43)			& 17	(22)		& -			& 27.0	&	-		&	-	\\
Slope					& 43	(32)		& 36 (38)			& 39 	(20)		& 65.5 (23.9)	& 44.9	& 43.8 (8.23)	& 44.8 (6.54)	\\
ANC$_0$					& 35	(11)		& 35	(23)			& 36 (13)		& 41.4 (7.56)	&31.6 	& 60.3 (3.27)	& 51.5 (3.61)	\\
\textbf{IOV
(CV\%)} 					& 			&				&			&			&		&			&	\\	
MMT	 					&	 -		& -				& 16	(8.5)		& -			& -		&	-		&	-	\\
\textbf{RUV} 				& 			&				&			&			&		&			&	\\	
exp.  (CV \%)				& 39.9		& 29.1			&			& 41.4		& 31.6	& 60.3 (3.27)	& 51.5 (3.61)	\\
add.  ($\cdot 10^9$ cells/L)	& 			& 0.626			&			&			& 		& 		& 	\\
box-cox.  					& 			& 				& 0.431	(2.6)	&			& 		& 		& 	\\
\\
\multicolumn{8}{l}{* baseline method was used, therefore, population median (and range) are given}\\
\multicolumn{8}{l}{$^\dagger$ the model also estimated the effect of Carboplatin on neutrophils}\\
\end{tabular}
\end{adjustbox}
\end{center}
\caption{Parameter estimates for different PD models describing paclitaxel-induced neutropenia. 
The model structure was introduced by Friberg et al (2002) \cite{Friberg2002}.
In subsequent publications the model structure (gold-standard) remained the same, but either the covariate model changed (Kloft et al. 2006 \cite{Kloft2006}), inter-occasion variability was modelled (Hansson et al. 2010 \cite{Hansson2010}) or it was fitted to a different patient population (Joerger et al. 2012 \cite{Joerger2012}, Henrich et al. 2017 \cite{Henrich2017}).
Later, the model was extended to also include the cumulative behavior of neutropenia (bone marrow exhaustion, BME) which could be observed over multiple cycles treatments (Henrich et al. 2017 \cite{Henrich2017}).
Note that in Kloft et al. (2006)  and Hansson et al. (2010), the Slope parameter was reported for unbound drug concentration (Slope = 54.5 and 69.6 L/$\mu$mol, respectively), which is assumed to be 5\% of the total drug concentration; thus we multiplied the values with 0.05 for the Slope parameter that is related to the total drug concentration (for a linear transformation, the relative standard error remains the same).
The typical values (TV) describe the fixed effects and the inter-individual variability (IIV) parameters the variability between patients. The residual variability (RUV) describes the deviation between measurements and model predictions accounting for measurement errors and potential model misspecifications. The IIV and RUV parameters are provided as coefficients of variation (CV). Note that baseline method B2 was used for baseline neutrophil counts $\text{ANC}_0$, i.e., the IIV was estimated together with RUV as one single parameter \cite{Dansirikul2008}.
NR: not reported}
\label{tab:PaclitaxelModelsALL}
\end{sidewaystable}

\subsection{Bone Marrow Exhaustion model}
In the CEPAC-TDM study \cite{Joerger2016}, cumulative neutropenia was observed, i.e., the lowest neutrophil concentration (nadir) as well as the maximum neutrophil concentration were decreasing over the course of treatment. 
A potential hypothesis for this cumulative behavior is that the drug also affects the long-term recovery of the bone marrow (bone marrow exhaustion). 
The structure of the gold-standard model for neutropenia by Friberg et al. \cite{Friberg2002} does not describe this long-term effect and was shown to overpredict neutrophil concentrations at later cycles \cite[section 3.3]{Henrich2017Thesis}. 
Therefore, Henrich et al. \cite{Henrich2017} have extended the model to include a stem cell compartment `Stem', representing pluripotent stem cells with slower proliferation which are also affected by the drug, see Figure~1 (blue part).


The proliferation rates for the two compartments `$\Prol$' and `$\Stem$' are given by
\begin{align*}
\kprol &= \ftr \cdot \ktr \\
k_\mathrm{stem} &= (1-\ftr) \cdot \ktr \,, \\
\end{align*}
respectively, where $\ktr = 4/\text{MTT}$. 
The baseline neutrophil count $\Circ_0 = \text{ANC}_0$ was inferred from the baseline data point (baseline method 2 \cite{Dansirikul2008})
\begin{equation*}
\Circ_{0,i}= y_{0} \cdot e^{\theta_\textrm{RV} \cdot \eta_{\Circ_0,i}}\,, \qquad \eta_{\Circ_0,i} \sim \mathcal{N}(0,1)\,.
\end{equation*}

The system of ODEs describing the structural model reads
 \begin{align*}
 \frac{d \Stem}{dt} &= k_{\text{stem}}\Stem \cdot (1-\Edrug) \cdot \biggl( \frac{\Circ_0}{\Circ} \biggr)^\gamma - k_\text{stem} \Stem \,, &\Stem(0) = \Circ_0 \\
\frac{d \Prol}{dt} &= k_{\text{prol}}\Prol \cdot (1-\Edrug) \cdot \biggr( \frac{\Circ_0}{\Circ} \biggl)^\gamma + k_\text{stem} \Stem - k_\text{tr} \Prol\,, &\Prol(0) = \Circ_0\\
\frac{d \text{Transit1}}{dt} &= k_\text{tr}\Prol - k_\text{tr} \text{Transit1}\,, &\text{Transit1}(0)=\Circ_0\\
\frac{d \text{Transit2}}{dt} &= k_\text{tr}\text{Transit1} - \ktr \text{Transit2}\,, &\text{Transit2}(0)=\Circ_0\\
\frac{d \text{Transit3}}{dt} &= k_\text{tr}\text{Transit2} - \ktr \text{Transit3}\,, &\text{Transit3}(0)=\Circ_0\\
\frac{d \Circ}{dt} &= \ktr \text{Transit3} - k_\text{circ} \Circ \,, &\text{Circ}(0)=\Circ_0\\
\end{align*}
with a linear drug effect, $\Edrug(t) = \text{Slope}\cdot C_1(t)$ and $k_\text{circ} = k_\text{tr}$.

Henrich et al. \cite{Henrich2017} estimated the model parameters in a population analysis based on the CEPAC-TDM study data \cite{Joerger2016}.

The typical model predictions in Figure~1 were generated for the typical patient of the CEPAC-TDM study: male, 56 years, $\text{ANC}_0 = 6.48 \cdot 10^9\, \text{cells}/\text{L}$, $\text{BSA}=1.8$, $\text{BILI}=7\, \mu \text{mol}/\text{L}$.

\section{Details on the continued learning approach}

\subsection{Derivation of the sequential Bayesian formulation}
\label{sec:deriv-full-conditionals}

The full hierarchical Bayesian model for individuals $1,...,N$ ($N_\text{TDM}$ in the main text) corresponds to the joint distribution of observations, individual parameters, and population parameters, denoted $p(y_{1:N},\theta_{1:N},\theta^\text{TV},\Omega)$, and which can be decomposed as
\begin{equation}
\label{eq:fullbayes}
p(y_{1:N},\theta_{1:N},\theta^\text{TV},\Omega) = p(\theta^\text{TV})p(\Omega)\prod_{i=1}^{N}\Big(p(y_{i}|\theta_{i}) p(\theta_{i}|\theta^\text{TV},\Omega) \Big),
\end{equation}
which can be written in the sequential formulation 
\begin{equation}
\label{eq:fullbayes-seq}
p(y_{1:i},\theta_{1:i},\theta^\text{TV},\Omega) = p(y_{i}|\theta_{i}) p(\theta_{i}|\theta^\text{TV},\Omega)  p(y_{1:i-1},\theta_{1:i-1},\theta^\text{TV},\Omega),\qquad i = 1,...,N.
\end{equation}
Then, the joint posterior $p(\theta_{i},\thetaTV,\Omega|y_{1:i})$ is given by
\begin{align*}
p(\theta_{i},\thetaTV,\Omega|y_{1:i}) &\propto \int p(\theta_{1:i},\thetaTV,\Omega,y_{1:i})d\theta_{1:i-1}\\
&= p(y_{i}|\theta_{i})p(\theta_{i}|\theta^\text{TV},\Omega) p(y_{1:i-1},\thetaTV,\Omega)\\
&\propto p(y_{i}|\theta_{i})p(\theta_{i}|\theta^\text{TV},\Omega) p(\thetaTV,\Omega|y_{1:i-1}).
\end{align*}
Hence, the joint posterior $p(\theta_{i},\thetaTV,\Omega|y_{1:i})$ in step $i$ depends on $y_{1:i-1}$ only through the marginal posterior $p(\thetaTV,\Omega|y_{1:i-1})$ from step $i-1$.
The conditional distributions in \Eqref{eq:condthetaTV},\Eqref{eq:condOmega} and \Eqref{eq:condTheta}, which are used in the Metropolis-Hastings-within-Gibbs sampling approach, are directly obtained from the above expression by dropping the constant terms in each case.


\subsection{Sampling in \Eqref{eq:condthetaTV} and \Eqref{eq:condOmega}}
\label{sec:CLapproach}

Sampling in iteration $l$ (Alg.~1; line 13) from \Eqref{eq:condthetaTV} corresponds (in our setting) to sampling from a multivariate normal distribution $\mathcal{N}(\mu^{\thetaTV(l)}_i,\Sigma^{\thetaTV(l)}_i)$ with parameters
\begin{align}\label{eq:cond_thetaTV_mean}
    \Sigma^{\thetaTV(l)}_i &= \biggl( \Big(\mathbf{S}^{\thetaTV}_{i-1}\Big)^{-1} + \Big(\Omega_i^{(l-1)}\Big)^{-1}\biggr)^{-1}\\
    \mu^{\thetaTV(l)}_i &= \Sigma^{\thetaTV}_i \biggl( \Big(\Omega_i^{(l-1)}\Big)^{-1} \theta_i^{(l-1)} + \Big(\mathbf{S}^{\thetaTV}_{i-1}\Big)^{-1}\overline{\thetaTV_{i-1}}\biggr)\,,
\end{align}
which corresponds to the update of a conjugate normal prior with a normal likelihood resulting in a normal posterior.
Equivalently, sampling from \Eqref{eq:condOmega} in Alg.~1; line 14 corresponds to sampling from an inverse-Wishart distribution $\mathcal{IW}(\Sigma^{\Omega(l)}_i,\nu^{(l)}_i)$ with parameters

\begin{align}
    \Sigma_i^{\Omega(l)} &= (\nu_i-d-1)\bar{\Omega}_{i-1} + \Big(\theta_i^{(l-1)}-\thetaTVl_i\Big)\Big(\theta_i^{(l-1)}-\thetaTVl_i\Big)^T \\
    \nu_i^{(l)} &= \nu_{i-1}+1\,.
    \label{eq:cond_omega2_nu}
\end{align}

\section{Additional analyses of the simulation study}
\label{sec:SupplementAnalyses}

\begin{figure}[H]
    \centering
    \includegraphics[width=\linewidth]{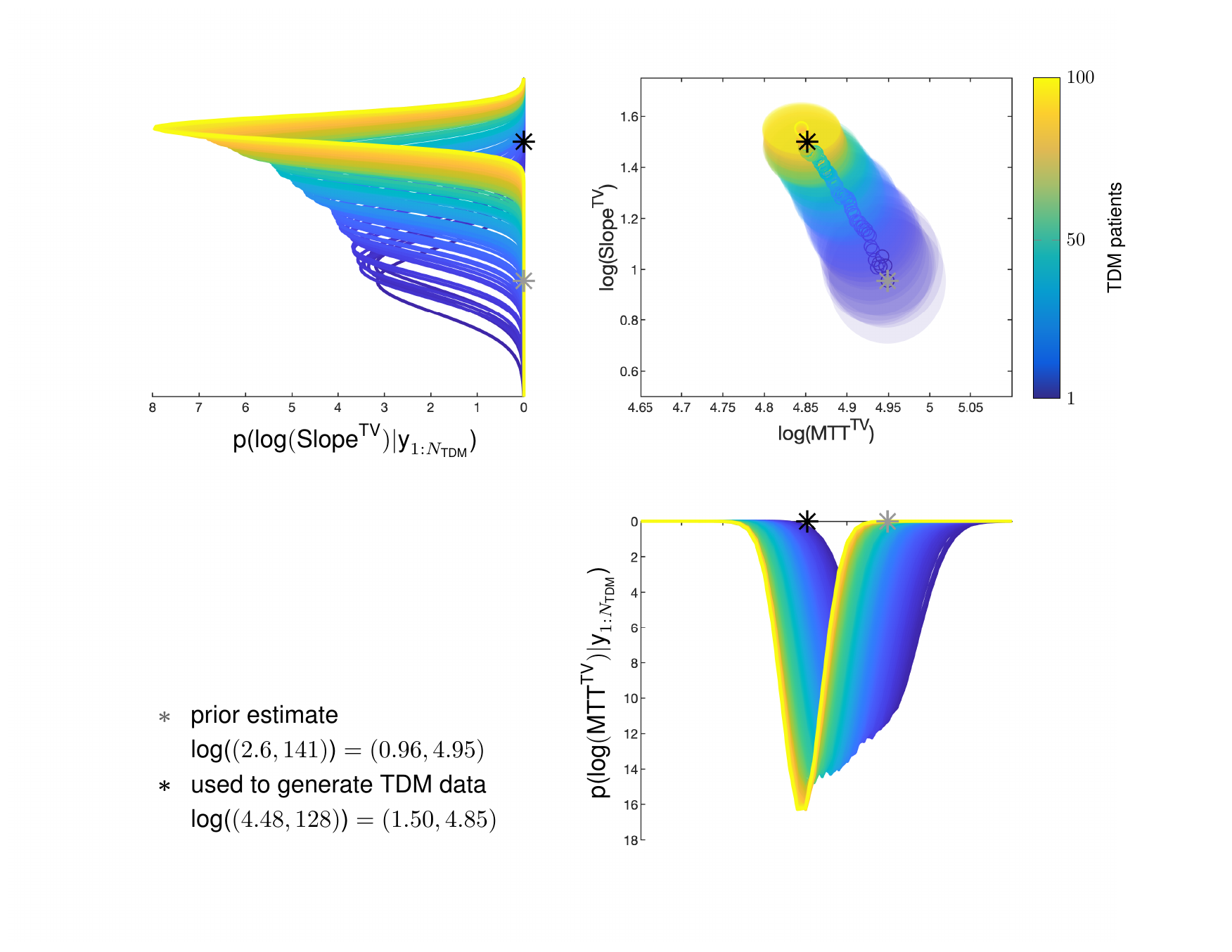}
    \caption{ Exemplary multivariate update of the typical values for MTT (mean transit time) and Slope (drug effect parameter) for the rich sampling scheme (TDM data every third day). The initial model used in the MIPD approach is the gold-standard model and the TDM data is generated using the gold-standard~R model. In this case also the parameter $\gamma$ was estimated on the individual level.}
    \label{fig:NAP_mvlogthetaTV1_100_gamma}
\end{figure}

\paragraph{Results including an estimation of $\gamma$.}
When $\gamma$ is included on the individual level inference, the estimation of the typical parameter values for `Slope' and `MTT' across patients is improved, see \ref{fig:NAP_mvlogthetaTV1_100_gamma}.

\begin{figure}[H]
\centering
\includegraphics[width=\textwidth]{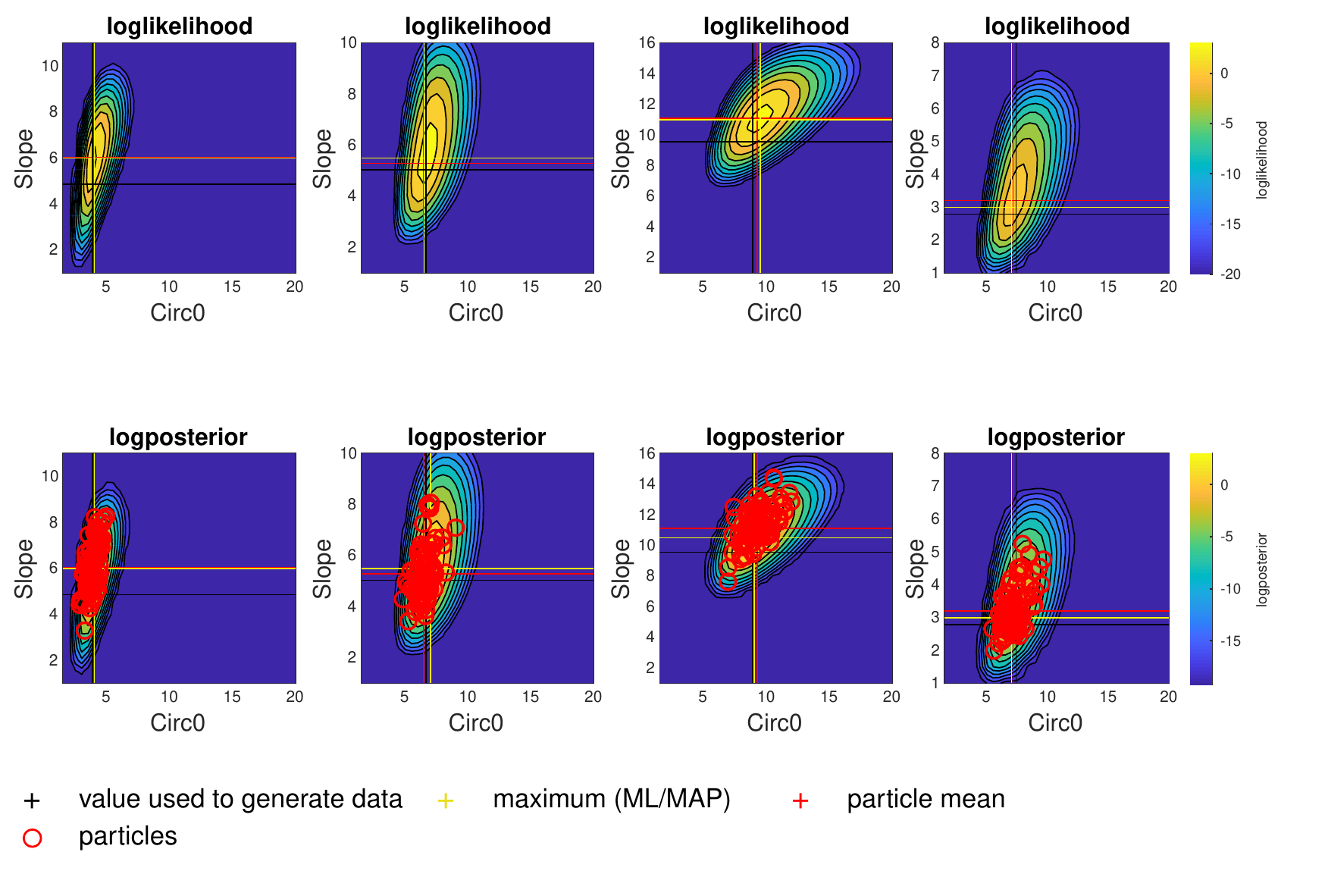}
\caption{Practical identifiability for the weekly sampling design. For illustration, $M'=100$ particles were resampled from the $M=1000$ particles used within the particle filter.}
\label{fig:Identifiability}
\end{figure}

\paragraph{Parameter identifiability.} 
To investigate the practical identifiability for weekly sampling (intermediate sampling scheme), we exemplarily computed the log-likelihood for four virtual patients at the end of the therapy, see \ref{fig:Identifiability}.
In order to exclude effects from other parameters we investigated a simplified setting in which only the parameters `Slope' and `Circ$_0$' were estimated on the individual level.
For some virtual patients, the log-likelihood takes the same values for various Slope values, which can be seen from the elongated yellow ranges covering a larger range of Slope values. 
In addition, the data suggest larger Slope values as the maximum of the likelihood (ML, large yellow cross) is reached for larger Slope values than used to generate the data (large black cross) in the upper panels. 

Furthermore, it can be observed that the prior only has minor influence on the log-posterior, when comparing the upper panels (log-likelihood) with the corresponding lower panels (log-posterior).
The analysis mean (large red cross) is close to the maximum-a-posteriori (MAP) estimate (large yellow cross in bottom panels).
Note that for skewed distributions the mean is different from the mode, and therefore, the particle mean should not be directly compared to the MAP.
The particles (red circles) cover areas of high posterior probability very well.

\begin{figure}[H]
\centering
\includegraphics[width=\textwidth]{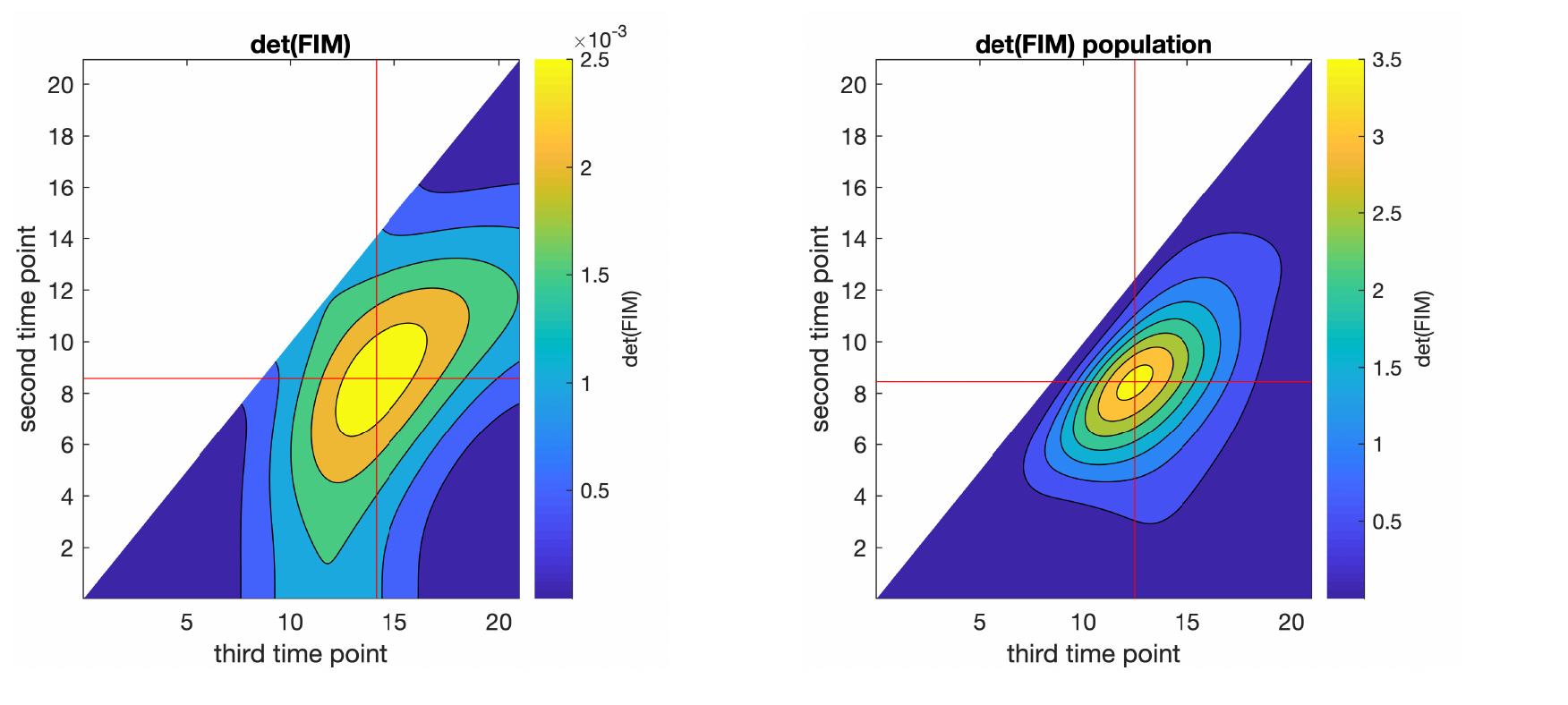}
\caption{Optimal design for the gold-standard model. The first sample timepoint at day 1 is fixed and only two more sample points within one cycle were allowed. The left panel shows the D-optimality criterion landscape for the typical patient and the right panel for the population.}
\label{fig:OptimalDesign}
\end{figure}

\paragraph{Optimal design.}
To investigate whether the reason for the practical identifiability are the chosen sampling time points, we investigated the optimal design for a design with three sampling timepoints where the first sampling time point at day 1 is fixed.
To infer the optimal design we used the frequently used criterion of D-optimality, i.e., choosing the design that maximizes the determinant of the Fisher information matrix.
The optimal design was determined for the typical patient (\ref{fig:OptimalDesign} left) and the whole patient population (\ref{fig:OptimalDesign} right).
The optimal second time point is approximately one day later than in the weekly sampling scheme (day 7) and when only the typical patient is considered, the third time point of the weekly sampling scheme (day 14) is chosen well. However, when the design is chosen based on the entire patient population, an earlier third time point is suggested.

\begin{figure}[H]
\centering
\includegraphics[width=\textwidth]{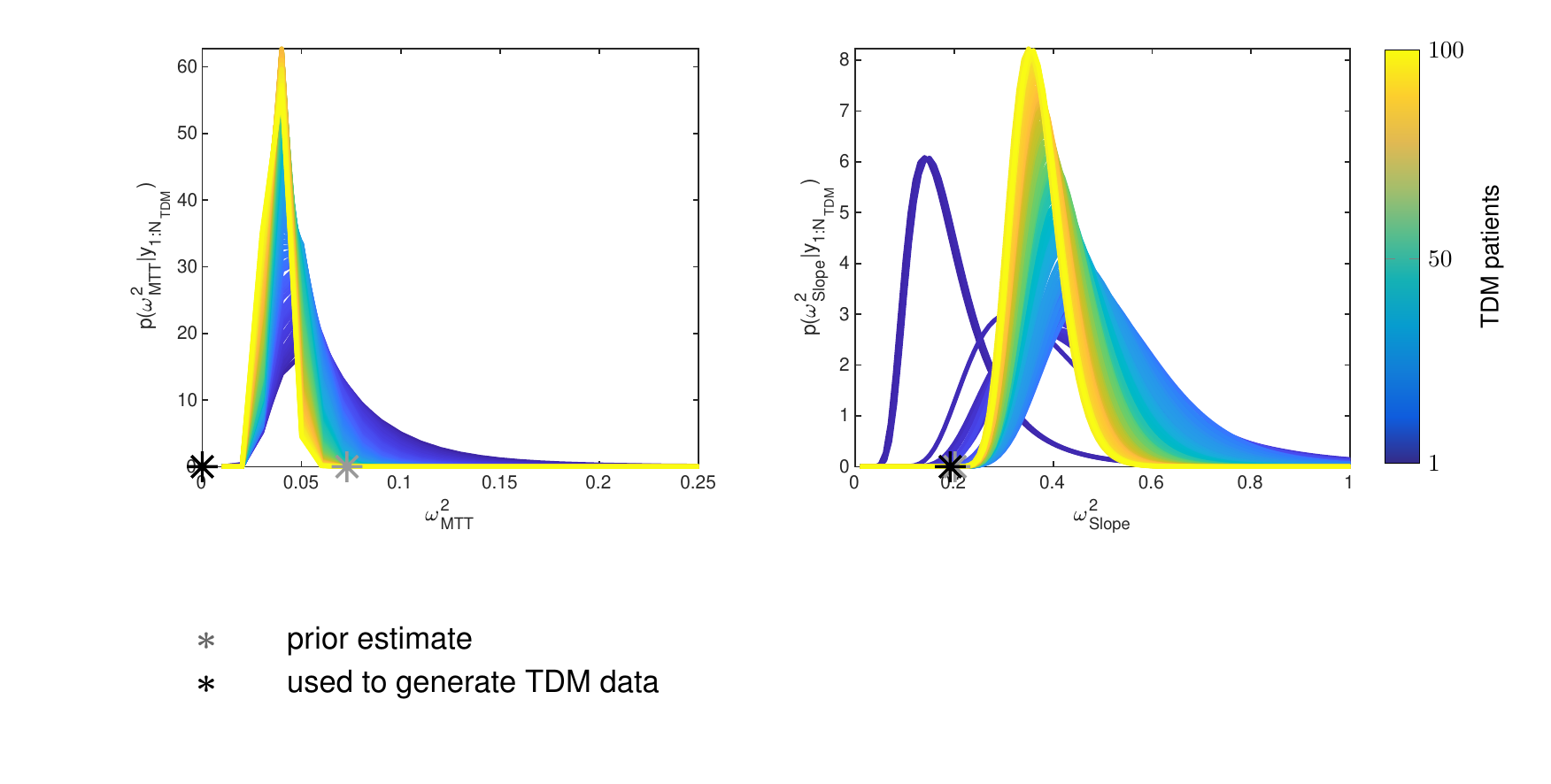}
\caption{Updates of the IIV parameters across patients.}
\label{fig:mvGMH_Omega}
\end{figure}

\paragraph{Updates of variability parameters.}
In addition to the typical values, also the variability parameters $\omega^2_\text{MTT}$ and $\omega^2_\text{Slope}$ were updated across patients, see \ref{fig:mvGMH_Omega}.
The IIV parameter for `MTT' moves from the prior estimate towards zero as no IIV has been estimated in the gold-standard~R model; in other words, TDM data was generated with the same parameter value for `MTT' for all patients.
The IIV parameter for `Slope' initially increased (dark blue) as individual estimated `Slope' parameters deviated considerably from the biased prior typical `Slope' parameter.
However, as more TDM data were observed and the typical value was increased, also the IIV parameter moved back towards the target value (close to the prior value).
A potential reason for the slight overestimation of the IIV parameters could be that we did not estimate the RUV parameter $\sigma$, but data were generated with an increased parameter value for $\sigma$. 

\paragraph{Parameter bias.}

\begin{figure}[H]
\includegraphics[width=\linewidth]{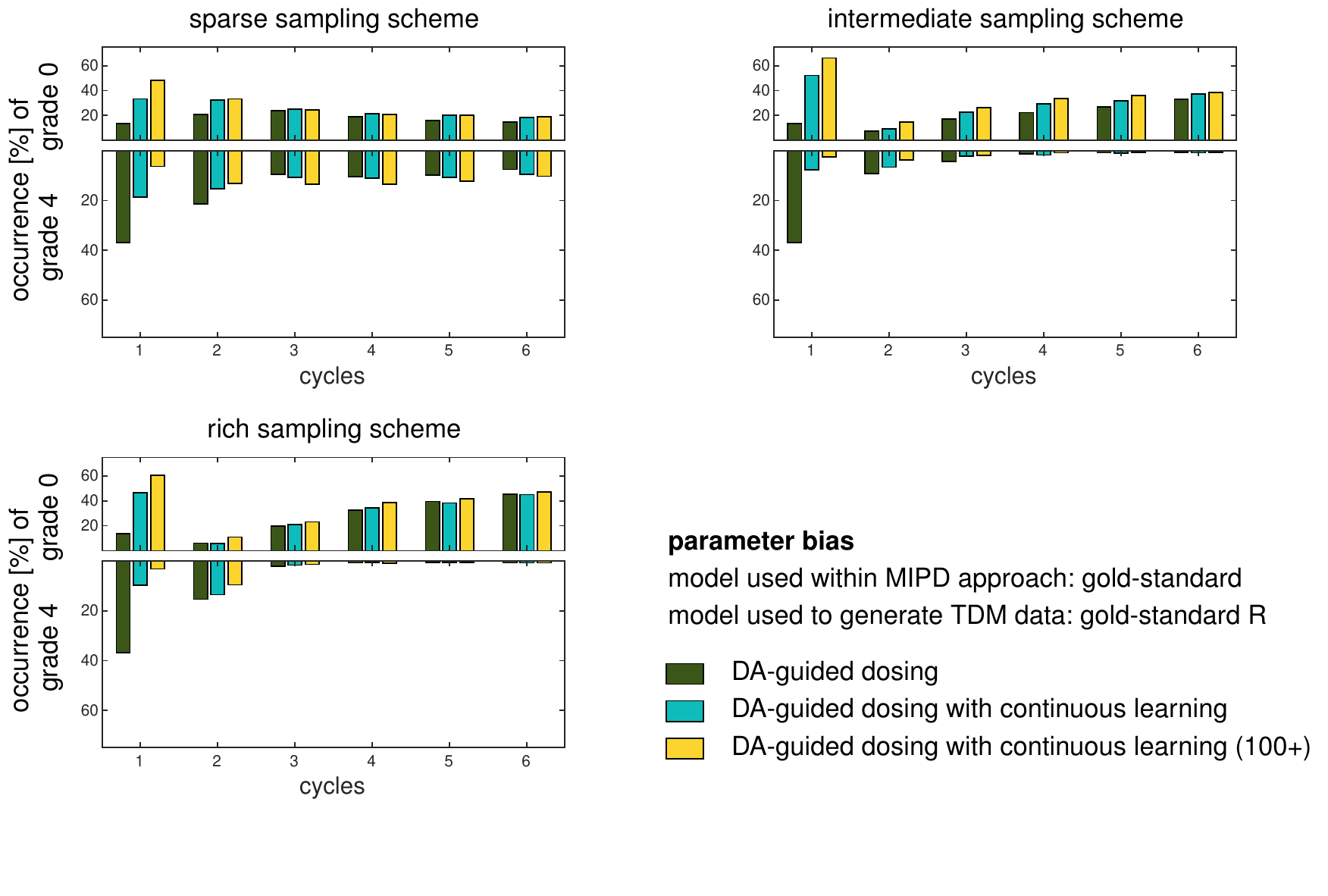}
\caption{Impact of continued learning approach on MIPD outcomes. TDM data was generated for $N_\text{TDM}=100$ virtual patients based on the gold-standard~R model. DA-guided dosing based on the gold-standard model was used alone or in conjunction with continued learning of the population parameters. The analysis was repeated $10$ times to account for statistical variability. 
}
\label{fig:ImpactonMIPD}
\end{figure}

Especially at the beginning of a patient's therapy when the misspecified prior dominated, MIPD benefited substantially from updating the model with every patient, see \ref{fig:ImpactonMIPD}. 
The occurrence of grade 4 in the first cycle was considerably reduced compared to DA-guided dosing alone.
As more patient-specific data were collected, individual parameters were increasingly well estimated with the DA approach and the influence of the misspecified prior vanished.
The occurrence of grade 0 is slightly increased in later cycles, which might be related to the overestimation of the variability parameters ($\omega^2_\text{MTT},\omega^2_\text{Slope}$), see \Figref{fig:mvGMH_Omega}.
We do not estimate the RUV parameter $\sigma$, however, the TDM data was generated with an increased parameter value for $\sigma$, therefore, the variability parameters capture to some extent the increased variability in the data.
This increases the uncertainty on the individual level.
DA-guided dosing alone has in this case the advantage that the IIV parameters are fairly similar between the gold-standard model and gold-standard~R model.


The population updates improved the used MIPD approach, especially for the first treatment cycle, when no patient-specific TDM data was available and the dose was solely determined based on {\it a priori} predictions.

\begin{figure}[H]
\centering
\includegraphics[width=\textwidth]{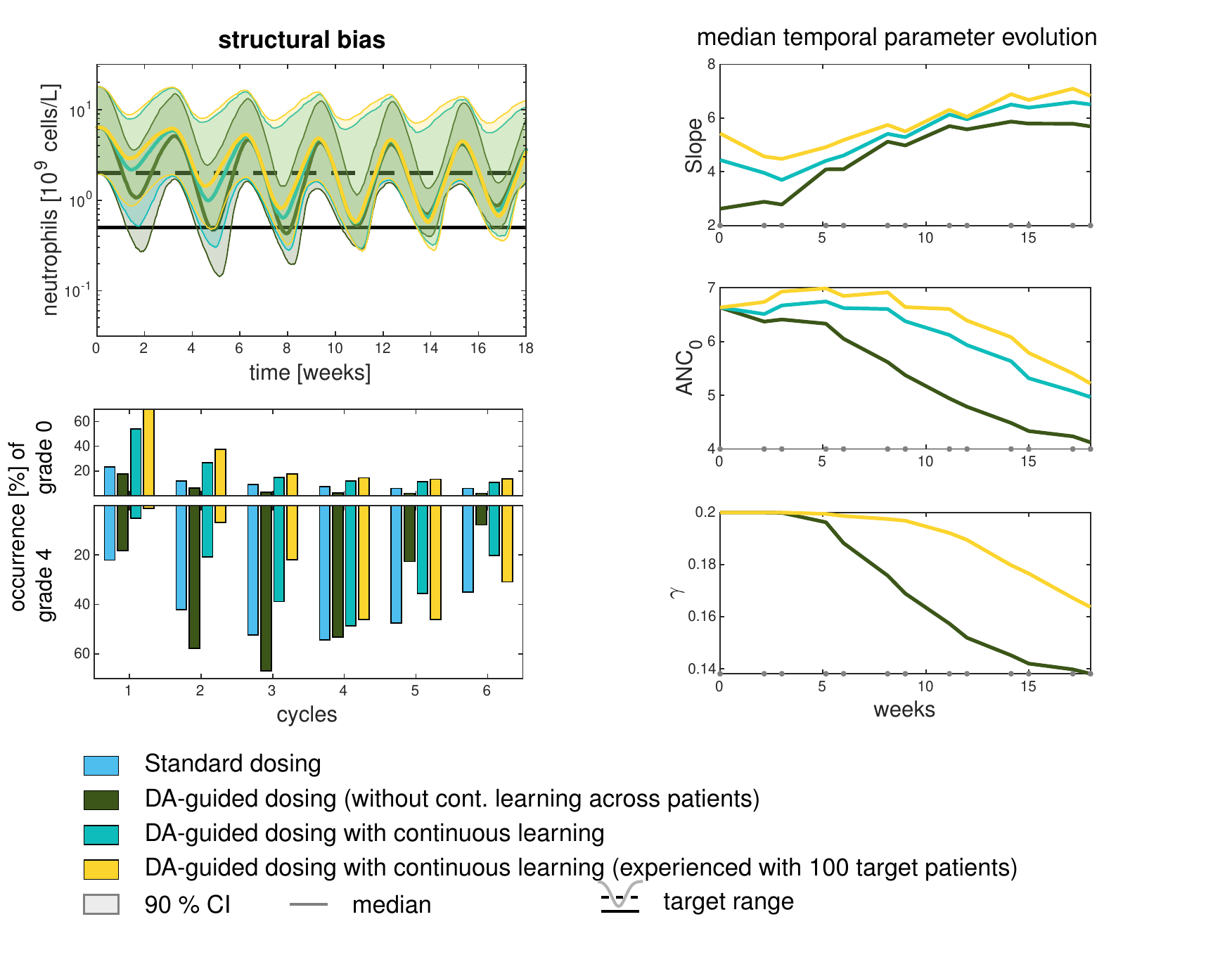}
\caption{Learning of temporal changes and adapting to structural model changes considering the \textbf{sparse sampling} scheme.}
\label{fig:MB3gammaSPARSE}
\end{figure}

\begin{figure}[H]
\centering
\includegraphics[width=\textwidth]{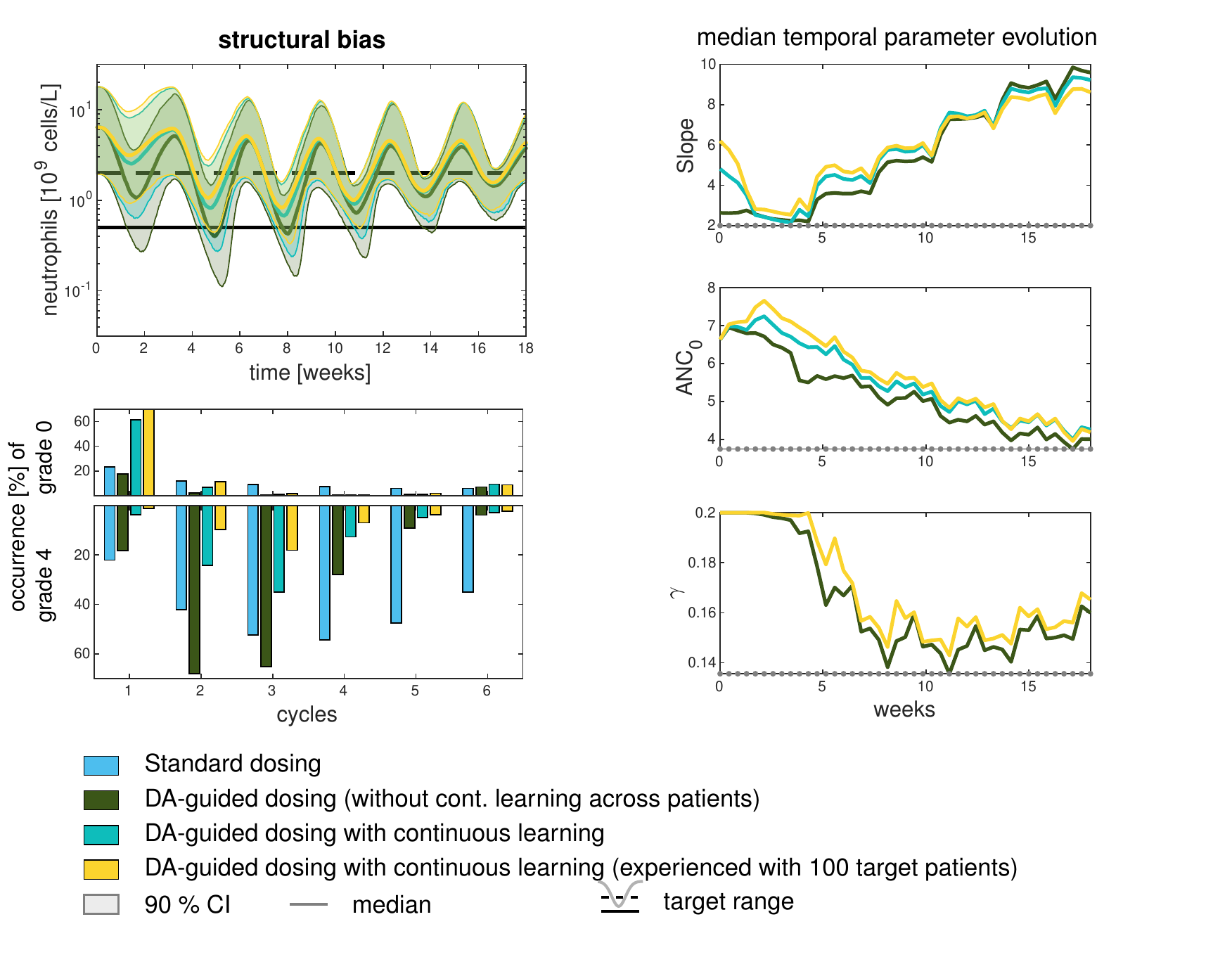}
\caption{Learning of temporal changes and adapting to structural model changes considering the \textbf{rich sampling} scheme.}
\label{fig:MB3gammaRICH}
\end{figure}

\paragraph{Learning of temporal changes for the sparse and rich sampling schemes.}
In the main text, learning of temporal changes was only shown for the intermediate (weekly) sampling scheme. For completeness, \ref{fig:MB3gammaSPARSE} and \ref{fig:MB3gammaRICH} show the results for the sparse and rich sampling schemes, respectively.
In the case of sparse TDM data, the continued learning updates decrease the incidence of grade 4 neutropenia in early cycles, but even lead to an increase in later cycles. This might be again related to the overestimation of the magnitude of inter-individual variability.
For the rich sampling scheme, the results are comparable to the intermediate sampling scheme presented in the main manuscript.
Note that a direct comparison of the parameter estimates is not recommended due to the structural differences of the models.

\end{document}